\documentclass[11pt]{amsart}

\usepackage{latexsym}
\usepackage{amssymb}
\usepackage{amsmath}

\newtheorem{theorem}{Theorem}[section]

\newtheorem{proposition}[theorem]{Proposition}
\newtheorem{corollary}[theorem]{Corollary}

\newcommand{\prend}{ $\diamondsuit $\hfill \bigskip}

\newcommand\supp{\mathop{\rm supp}}

\newcommand\esssup{\mathop{\rm esssup}}

\newcommand\id{\mathop{\rm id}}

\newcommand\rank{\mathop{\rm rank}}

\newcommand\vn{\mathop{\rm VN}}
\newcommand\nph{\varphi}
\newcommand\dist{\mathop{\rm dist}}

\newcommand\hh{{\mathop{\rm h}}}
\newcommand{\haag}{\otimes_\hh}
\newcommand{\haags}{\haag\cdots\haag}
\newcommand\op{\mathop{\rm op}}
\newcommand\dd{\mathop{\rm d}}
\newcommand\mm{\mathop{\rm m}}
\newcommand\eh{{\mathop{\rm eh}}}
\newcommand\ph{\mathop{\rm ph}}
\newcommand{\A}{\mathcal{A}}
\newcommand{\As}[1][\A]{{#1}_1,\dots,{#1}_n}
\newcommand{\timess}{\times\cdots\times}
\newcommand{\otimess}{\otimes\cdots\otimes}
\newcommand{\cl}[1]{\mathcal{#1}}
\newcommand{\bb}[1]{\mathbb{#1}}
\newcommand{\ehaag}{\otimes_\eh}
\newcommand{\ehaags}{\ehaag\cdots\ehaag}

\newcommand{\defeq}{\stackrel{\mathop{\rm def}}=}
\newcommand{\odots}{\odot\cdots\odot}

\newcommand{\B}{\mathcal{B}}

\newcommand{\E}{\mathcal{E}}
\newcommand{\F}{\mathcal{F}}

\newcommand{\K}{\mathcal{K}}

\newcommand{\X}{\mathcal{X}}
\newcommand{\Y}{\mathcal{Y}}
\newcommand{\CC}{\mathit{CC}}

\newcommand{\ff}{\mathit{f\!f}}
\newcommand{\Kh}{\K_\hh}
\begin{document}

\title{Schur and operator multipliers}

\author{I. G. Todorov and L. Turowska}

\date{}
\thanks{The first author was supported by EPSRC grant D050677/1.
The second author was supported by the Swedish Research Council.}

\maketitle

Schur multipliers were introduced by Schur in the early 20th century
and have since then found a considerable number of applications in
Analysis and enjoyed an intensive development. Apart from the beauty
of the subject in itself, sources of interest in them were
connections with Perturbation Theory, Harmonic Analysis, the Theory
of Operator Integrals and others. Advances in the quantisation of
Schur multipliers were recently made in \cite{ks}. The aim of the
present article is to summarise a part of the ideas and results
in the theory of Schur and
operator multipliers. We start with the classical Schur multipliers
defined by Schur and their characterisation by Grothendieck, and
make our way through measurable multipliers studied by Peller and
Spronk, operator multipliers defined
by Kissin and Shulman and, finally, multidimensional Schur and
operator multipliers developed by Juschenko and the authors. We
point out connections of the area with Harmonic Analysis and the
Theory of Operator Integrals.

\section{Classical Schur multipliers}

For a Hilbert space $H$, let $\cl B(H)$ be the collection of all bounded linear
operators acting on $H$ equipped with its operator norm $\|\cdot\|_{\op}$.
We denote by $\ell^2$ the Hilbert space of all square summable
complex sequences. With an operator $A\in \cl B(\ell^2)$, one can associate a matrix
$(a_{i,j})_{i,j\in \bb{N}}$ by letting $a_{i,j} = (Ae_j,e_i)$, where $\{e_i\}_{i\in \bb{N}}$
is the standard orthonormal basis of $\ell^2$. The space $M_{\infty}$ of all
matrices obtained in this way is a subspace of the space $M_{\bb{N}}$ of all
complex matrices indexed by $\bb{N}\times\bb{N}$. It is easy to see that the correspondence between
$\cl B(\ell^2)$ and $M_{\infty}$ is one-to-one.

Any function $\nph : \bb{N}\times\bb{N}\rightarrow\bb{C}$ gives rise to a linear transformation
$S_{\nph}$ acting on $M_{\bb{N}}$ and given by $S_{\nph}((a_{i,j})_{i,j}) = (\nph(i,j)a_{i,j})_{i,j}$.
In other words, $S_{\nph}((a_{i,j})_{i,j})$ is the entry-wise product of the matrices
$(\nph(i,j))_{i,j}$ and $(a_{i,j})_{i,j}$, often called {\bf Schur product}.
The function $\nph$ is called a {\bf Schur multiplier} if $S_{\nph}$ leaves the subspace $M_{\infty}$
invariant.
We denote by $S(\bb{N},\bb{N})$ the set of all Schur multipliers.

Let $\nph$ be a Schur multiplier. Then the correspondence between
$\cl B(\ell^2)$ and $M_{\infty}$ gives rise to a mapping (which we
denote in the same way) on $\cl B(\ell^2)$. We first note that
$S_\nph$ is necessarily bounded in the operator norm.
This follows from the Closed Graph Theorem;
indeed, suppose that $A_k\rightarrow 0$ and
$S_{\nph}(A_k)\rightarrow B$ in the operator norm, for some elements
$A_k, B\in \cl B(\ell^2)$, $k\in \bb{N}$. Letting $(a^k_{i,j})$ and
$(b_{i,j})$ be the corresponding matrices of $A_k$ and $B$, we have
that $a_{i,j}^k = (A_ke_j,e_i)\rightarrow_{k\rightarrow \infty} 0$,
for each $i,j\in \bb{N}$. But then $\{\nph(i,j)a^k_{i,j}\}_{k\in
\bb{N}}$ converges to both $b_{i,j}$ and $0$, and since this holds
for every $i,j\in \bb{N}$, we conclude that $B = 0$. Let
$\|\nph\|_{\mm}$ denote the norm of $S_{\nph}$ as a bounded operator
on $\cl B(\ell^2)$; we call $\|\nph\|_{\mm}$ the multiplier norm of
$\nph$.

It now follows that $S(\bb{N},\bb{N})\subseteq \ell^{\infty}(\bb{N}\times\bb{N})$.
Indeed, if $\{E_{i,j}\}_{i,j\in
\bb{N}}$ is the canonical matrix unit system in $\cl B(\ell^2)$ then
we have that $|\nph(i,j)| = \|S_{\nph}(E_{i,j})\|_{\op}\leq
\|\nph\|_{\mm}$, for all $i,j\in \bb{N}$.
It is trivial to verify that
$S(\bb{N},\bb{N})$ is a subalgebra of $\ell^{\infty}(\bb{N}\times\bb{N})$ when the latter is equipped
with the usual pointwise operations;
moreover, it can be shown that $S(\bb{N},\bb{N})$ is
a semi-simple commutative Banach algebra when equipped with the norm $\|\cdot\|_{\mm}$.
We note in passing that, since $\nph$ is a bounded function, the
restriction of $S_{\nph}$ to the class $\cl C_2(\ell^2)$ of all
Hilbert-Schmidt operators is bounded when $\cl C_2(\ell^2)$ is
equipped with its Hilbert-Schmidt norm $\|\cdot\|_2$, and its norm
as an operator on $\cl C_2(\ell^2)$ is equal to $\|\nph\|_{\infty}$.
This follows from the fact that if $A\in \cl C_2(\ell^2)$ and
$(a_{i,j})_{i,j}\in M_{\infty}$ is the matrix corresponding to $A$
then $\|A\|_2 = \left(\sum_{i,j} |a_{i,j}|^2\right)^{\frac{1}{2}}$.

Another immediate observation is the fact that if
$(a_{i,j})_{i,j}\in M_{\infty}$ and $a\in
\ell^{\infty}(\bb{N}\times\bb{N})$ is the function given by $a(i,j)
= a_{i,j}$, then $a$ is a Schur multiplier and $\|a\|_{\mm}\leq
\|a\|_{\op}$. To see this, let $b\in M_{\infty}$ and $B\in \cl
B(\ell^2)$ be the operator corresponding to $b$. If $A\in \cl
B(\ell^2)$ is the operator corresponding to $a$, we have by general
operator theory that the norm of the operator $A\otimes B \in \cl
B(\ell^2\otimes\ell^2)$ is equal to $\|A\|_{\op}\|B\|_{\op}$.
Identifying $\ell^2\otimes\ell^2$ with $\ell^2(\bb{N}\times\bb{N})$,
and letting $P\in \cl B(\ell^2\otimes\ell^2)$ be the projection on
the closed linear span of $\{e_i\otimes e_i\}_{i\in \bb{N}}$, we see
that the matrix $(a_{i,j}b_{i,j})_{i,j}$ corresponds to the operator
$P(A\otimes B)P$. Thus, $\|S_a(B)\|_{\op} = \|P(A\otimes
B)P\|_{\op}\leq \|A\|_{\op}\|B\|_{\op}$, and the claim follows. In
fact, $M_{\infty}$, equipped with the multiplier norm
$\|\cdot\|_{\mm}$ is a (commutative) semi-simple Banach subalgebra
of $S(\bb{N},\bb{N})$.

We summarise the above inclusions:
$$M_{\infty}\subseteq S(\bb{N},\bb{N})\subseteq \ell^{\infty}(\bb{N}\times\bb{N}).$$
Both of them are strict: for the first one, take the constant
function ${\bf 1}$ taking the value $1$ on $\bb{N}\times\bb{N}$.
Obviously, ${\bf 1}$ is a Schur multiplier (in fact, $S_{{\bf 1}}$
is the identity transformation) but does not belong to $M_{\infty}$
since the rows and columns of the matrices in $M_{\infty}$ are
square summable. The fact that the second inclusion is strict is
much more subtle. An example of a function which belongs to
$\ell^{\infty}(\bb{N}\times\bb{N})$ but not to $S(\bb{N},\bb{N})$ is
the characteristic function $\chi_{\Delta}$ of the set $\Delta =
\{(i,j) : j\leq i\}$, see for example \cite{dav}. 

The question that arises is: \lq\lq which'' functions are Schur
multipliers? The following description was obtained by Grothendieck
\cite{Gro}:

\begin{theorem}\label{groth}
Let $\nph\in \ell^{\infty}(\bb{N}\times\bb{N})$. The following are equivalent:

(i) \ $\nph$ is a Schur multiplier and $\|\nph\|_{\mm} < C$;

(ii) There exist families $\{a_k\}_{k\in \bb{N}}, \{b_k\}_{k\in \bb{N}}\in \ell^{\infty}$
such that
$$\sup_{i\in \bb{N}}\sum_{k=1}^{\infty} |a_k(i)|^2 < C, \ \sup_{j\in \bb{N}}\sum_{k=1}^{\infty} |b_k(j)|^2 < C$$
and
$$\nph(i,j) = \sum_{k=1}^{\infty} a_k(i)b_k(j), \ \ \ \ \mbox{ for all } i,j\in \bb{N}.$$
\end{theorem}

Suppose that we are given two {\it finite} families
$\{a_k\}_{k=1}^{N}, \{b_k\}_{k=1}^{N}$ $\subseteq$ $\ell^{\infty}$
and let $\nph\in \ell^{\infty}(\bb{N}\times\bb{N})$ be the function
given by $\nph(i,j) = \sum_{k=1}^{N} a_k(i)b_k(j)$. For $a\in
\ell^{\infty}$, let $D_a\in \cl B(\ell^2)$ be the operator whose
matrix has the sequence $a$ down its main diagonal and zeros
everywhere else. An easy computation shows that $\nph$ is a Schur
multiplier and that, in fact,
$$S_{\nph}(T) = \sum_{k=1}^N D_{a_k}TD_{b_k}, \ \ \ T\in \cl B(\ell^2).$$
The transformations on $\cl B(\ell^2)$ obtained in this way belong
to the important class of {\it elementary operators}. The norm of
this operator, and hence $\|\nph\|_{\mm}$, is bounded by
$(\sup_{i\in \bb{N}}\sum_{k=1}^{N} |a_k(i)|^2 \sup_{j\in
\bb{N}}\sum_{k=1}^{N} |b_k(j)|^2)^{1/2}$. In fact, for $f$, $g\in
\ell^2$
\begin{eqnarray}\label{estim}
|(S_{\nph}(T)f,g)| &=& \left|\left(\sum_{k=1}^N
D_{a_k}TD_{b_k}f,g\right)\right| \leq \sum_{k=1}^N |(TD_{b_k}f,D_{a_k}^*g)|\\
& \leq & \sum_{k=1}^N\|TD_{b_k}f\|\|D_{a_k}^*g\|\nonumber\\
& \leq & \left(\sum_{k=1}^N\|T\|^2\|D_{b_k}f\|^2\right)^{1/2} \left(\sum_{k=1}^N\|D_{a_k}^*g\|^2\right)^{1/2}\nonumber\\
& = &\|T\|\left(\sum_{k=1}^ND_{b_k}^*D_{b_k}f,f\right)^{1/2}\left(\sum_{k=1}^N D_{a_k}D_{a_k}^*g,g\right)^{1/2}\nonumber\\
& \leq &\|T\|\left\|\sum_{k=1}^{N}
D_{b_k}^*D_{b_k}\right\|^{1/2}\left\|\sum_{k=1}^{N}
D_{a_k}D_{a_k}^*\right\|^{1/2}\|f\|\|g\|\nonumber
\end{eqnarray}
which proves the claim.

If $\nph$ is an arbitrary Schur multiplier and
$\{a_k\}_{k\in \bb{N}}$ and $\{b_k\}_{k\in \bb{N}}$ are the families
arising in Grothendieck's description, then letting $\nph_N$ be the
function given by $\nph_N(i,j) = \sum_{k=1}^{N} a_k(i)b_k(j)$, $N\in
\bb{N}$, we have that $\nph_N\rightarrow \nph$ pointwise and,
moreover (by Grothendieck's theorem again), the norms
$\|\nph_N\|_{\mm}$ are uniformly bounded in $N$. Thus,
Grothendieck's characterisation can be viewed as a uniform
approximation result for Schur multipliers by \lq\lq elementary
Schur multipliers''.

There is another convenient formulation of Theorem \ref{groth} which
links the subject of Schur multipliers to Operator Space Theory. Namely,
the space of Schur multipliers can be identified with the extended Haagerup tensor product
$\ell^{\infty}\otimes_{\eh}\ell^{\infty}$. For the time being, let us take this statement as
the definition of the space $\ell^{\infty}\otimes_{\eh}\ell^{\infty}$; more on this
will be said later.

Truncation of matrices has been important in applications. Suppose that we are given
a matrix --- either of finite size, or an element of $M_{\infty}$. To truncate it
along a certain subset $\kappa\subseteq \bb{N}\times \bb{N}$ means
to replace it by the matrix that has the same entries on the subset $\kappa$
and zeros everywhere else.
Obviously, this is precisely the operation of
Schur-multiplying a matrix $a\in M_{\infty}$ by $\chi_{\kappa}$.
Thus, truncation along $\kappa$ is a well-defined (and automatically bounded)
transformation on $\cl B(\ell^2)$ if and only if $\chi_{\kappa}$ is a Schur multiplier.
The Schur multipliers of this form are precisely the idempotent elements of
the Banach algebra $S(\bb{N},\bb{N})$.

It is easy to exhibit idempotent Schur multipliers: if $\{\alpha_k\}_{k=1}^{\infty}$ and
$\{\beta_k\}_{k=1}^{\infty}$ are families of pairwise disjoint subsets of $\bb{N}$ and
$\kappa = \cup_{k=1}^{\infty} \alpha_k\times\beta_k$ then
$$S_{\chi_{\kappa}}(T) = \sum_{k=1}^{\infty} P_kTQ_k, \ \ \ T\in \cl B(\ell^2),$$
where $P_k$ (resp. $Q_k$) is the projection onto the closed span of
$\{e_i : i\in \alpha_k\}$ (resp. $\{e_i : i\in \beta_k\}$), and
hence $\|\chi_{\kappa}\|_{\mm}\leq 1$. Since $S(\bb{N},\bb{N})$ is
an algebra with respect to pointwise multiplication, the function
$\chi_{\kappa}$ is a Schur multiplier for every set $\kappa$
belonging to the subset ring generated by the sets of the above
form.
It is natural to ask whether the elements of this algebra are all
idempotent Schur multipliers. Although the answer to this question
is not known (in fact, the question is one of the difficult open
problems in the area), a related result was recently established by
Davidson and Donsig \cite{dd}. Let us call a subset $\kappa\subseteq
\bb{N}\times\bb{N}$ a {\bf Schur bounded pattern} if every function
$\nph\in \ell^{\infty}(\bb{N}\times\bb{N})$ supported on $\kappa$ is
a Schur multiplier. Obviously, if $\kappa$ is a Schur bounded
pattern then $\chi_{\kappa}$ is a Schur multiplier, and not vice
versa (just take $\kappa = \bb{N}\times\bb{N}$). The aforementioned
theorem reads as follows:

\begin{theorem}\label{davdon}
Let $\kappa\subseteq \bb{N}\times\bb{N}$. The following are equivalent:

(i) \ $\kappa$ is a Schur bounded pattern;

(ii) there exist sets $\kappa_r,\kappa_c\subseteq
\bb{N}\times\bb{N}$ and a number $N\in \bb{N}$ such that $\kappa_r$
(resp. $\kappa_c$) has at most $N$ entries in every row (resp.
column) and $\kappa = \kappa_r \cup \kappa_c$.
\end{theorem}

In fact, the Davidson and Donsig results cover the more general case of 
non-negative functions $\nph\in \ell^{\infty}(\bb{N}\times\bb{N})$
which are \lq\lq hereditary Schur multipliers'' (the terminology is ours) in the sense that
if $\psi\in \ell^{\infty}(\bb{N}\times\bb{N})$ and $|\psi|\leq \nph$ then $\psi$ is a Schur multiplier.

\bigskip

Once there is a complete characterisation of Schur multipliers, it
is natural to ask for a description of certain special classes of
Schur multipliers. Say that a Schur multiplier $\nph$ is compact if
the mapping $S_{\nph} : \cl B(\ell^2)\rightarrow \cl B(\ell^2)$ is a
compact operator. The following result was established by Hladnik in
\cite{hladnik}:

\begin{theorem}\label{hladnik}
A Schur multiplier $\nph\in \ell^{\infty}(\bb{N}\times\bb{N})$ is
compact if and only if a representation can be chosen for $\nph$ as
in Theorem \ref{groth} (ii), such that $\{a_k\}_{k\in \bb{N}},$
$\{b_k\}_{k\in \bb{N}}\subseteq c_0$ and the series
$\sum_{k=1}^{\infty} |a_k|^2$ and $\sum_{k=1}^{\infty} |b_k|^2$ are
uniformly convergent.
\end{theorem}

Positivity is another natural property that a Schur multiplier may
(or may not) have: a Schur multiplier $\nph$ is called positive if
$S_{\nph}(A)$ is a positive operator for every positive operator
$A$. The following holds:

\begin{theorem}
A Schur multiplier $\nph\in \ell^{\infty}(\bb{N}\times\bb{N})$ is positive
if and only if there exists
$\{a_k\}_{k\in \bb{N}}\subseteq \ell^{\infty}$ such that
$\sup_{i\in \bb{N}}\sum_{k=1}^{\infty} |a_k(i)|^2 < \infty$ and $\nph(i,j) = \sum_{k=1}^{\infty} a_k(i)\overline{a_k(j)}$,
for all $i,j\in \bb{N}$.
\end{theorem}

\section{Schur multipliers over measure spaces}

If $\nph\in \ell^{\infty}(\bb{N}\times\bb{N})$
then the operator on the Hilbert space $\ell^2(\bb{N}\times\bb{N})$
of multiplication by $\nph$ is bounded. If we equip $\ell^2(\bb{N}\times\bb{N})$ with the norm
arising from its identification with the space $\cl C_2(\ell^2)$ of all Hilbert-Schmidt operators
on $\ell^2$, then it is easy to see that $\nph$ is a Schur multiplier if and only if
this multiplication operator is bounded in the
operator norm. This approach is useful because
it allows us to study Schur multipliers in a more general setting.
To describe this setting, let $(X,\mu)$ and $(Y,\nu)$ be standard $\sigma$-finite measure spaces.
We equip $X\times Y$ with the product measure $\mu\times\nu$.
The space $L^2(X\times Y)$ can be canonically identified with the space $\cl C_2(L^2(X),L^2(Y))$
of all Hilbert-Schmidt operators from $L^2(X)$ into $L^2(Y)$: if $f\in L^2(X\times Y)$,
let $T_f$ be the Hilbert-Schmidt operator given by
$$T_f\xi (y) = \int_X f(x,y)\xi(x) d\mu(x), \ \ \ \xi\in L^2(X).$$
For $f\in L^2(X\times Y)$, let $\|f\|_{\op}$ be the operator norm of $T_f$.

Now let $\nph\in L^{\infty}(X\times Y)$. The operator $S_\nph :
L^2(X\times Y)\rightarrow L^2(X\times Y)$ of multiplication by
$\nph$ is bounded in the $L^2$-norm (its norm is equal to
$\|\nph\|_{\infty}$). If $S_{\nph}$ is moreover bounded in
$\|\cdot\|_{\op}$, that is, if there exists $C > 0$ such that
$\|\nph f\|_{\op}\leq C\|f\|_{\op}$ for every $f\in L^2(X\times Y)$,
then we call $\nph$ a Schur $\mu,\nu$-multiplier (or simply a Schur
multiplier if the measures are clear from the context). If $X = Y =
\bb{N}$ is equipped with the counting measure, this new notion
reduces to the one described in the previous section.

We note that the property of a function $\nph$ to be or not to be
a Schur multiplier depends only on the values of the function up to a {\it null
with respect to the product measure} set.

If $\nph\in L^{\infty}(X\times Y)$ is a Schur multiplier then the mapping $S_{\nph}$
extends by continuity to a mapping on $\cl K(L^2(X),L^2(Y))$; after taking
its second dual we arrive at a bounded weak* continuous
linear transformation (which we denote in the same
way) on $\cl B(L^2(X),L^2(Y))$. The multiplier norm $\|\nph\|_{\mm}$ of $\nph$ is defined as
the operator norm of $S_{\nph}$. We denote by $S(X,Y)$ the set of all Schur multipliers.
Clearly, $S(X,Y)$ is a subalgebra of $L^{\infty}(X\times Y)$ (with respect to
pointwise multiplication).

Let $\cl D_1$ (resp. $\cl D_2$) be the multiplication masa of $L^{\infty}(X)$
(resp. $L^{\infty}(Y)$). We denote by $M_a$ the element of $\cl D_1$ corresponding to
the element $a\in L^{\infty}(X)$; we use a similar notation for the operators in $\cl D_2$.
If $f\in L^2(X\times Y)$, $a\in L^{\infty}(X)$ and $b\in L^{\infty}(Y)$, then
$M_bT_fM_a = T_{f(a\otimes b)}$, where $(a\otimes b) (x,y) = a(x)b(y)$, $x\in X$, $y\in Y$.
It follows that $S_{\nph}$ is a $\cl D_2,\cl D_1$-bimodule map in the sense that
$S_{\nph}(BTA) = BS_{\nph}(T)A$, for all $T\in \cl B(L^2(X),L^2(Y))$, $A\in \cl D_1$
and $B\in \cl D_2$.

A characterisation similar to the Grothendieck's one also holds in
the measurable setting: the following result was established
by Peller \cite{peller_two_dim} (see also Spronk \cite{spronk}).

\begin{theorem}\label{groth2}
Let $\nph\in L^{\infty}(X\times Y)$. The following are equivalent:

(i) \ $\nph$ is a Schur multiplier and $\|\nph\|_{\mm} < C$;

(ii) there exist families $\{a_k\}_{k\in \bb{N}}\subseteq
L^{\infty}(X)$ and $\{b_k\}_{k\in \bb{N}}\subseteq L^{\infty}(Y)$
such that
$$\esssup_{x\in X}\sum_{k=1}^{\infty} |a_k(x)|^2 < C, \ \esssup_{y\in Y}\sum_{k=1}^{\infty} |b_k(y)|^2 < C$$
and
$$\nph(x,y) = \sum_{k=1}^{\infty} a_k(x)b_k(y), \ \ \ \mbox{for allmost all } (x,y)\in X\times Y.$$
\end{theorem}

We outline the proof of the theorem using  results of Smith
\cite{smith} and Haagerup \cite{haagerup}; this proof relies on the
notion of complete boundedness whose mention is deliberately omitted
for the time being but will be discussed in detail later.

\medskip

{\noindent {\it {Sketch of proof. }}} Assume $(X,\mu)=(Y,\nu)$ and
$\nph$ is a Schur $\mu,\mu$-multiplier with $\|\nph\|_{\mm}< C$.
Then $S_\nph$ can be extended to a bounded $\cl D,\cl D$-bimodule
map on $\cl K(L^2(X))$, where $\cl D$ is the multiplication masa of
$L^{\infty}(X)$. By \cite[Theorem 2.1, Theorem 3.1]{smith}, there
exist sequences $\{a_k\}_{k=1}^{\infty}$ and
$\{b_k\}_{k=1}^{\infty}$ in $L^{\infty}(X)$ such that
$\|\sum_{k=1}^{\infty} M_{b_k}M_{b_k}^*\|\|\sum
M_{a_k}^*M_{a_k}\|=\|\nph\|_{\mm}^2$ and, for all $T\in \cl
K(L^2(X))$,
\begin{equation}\label{sum}
S_\nph(T)=\sum_{k=1}^{\infty} M_{b_k}TM_{a_k},
\end{equation}
where the sum on the
right hand side converges in the strong operator topology.
Let $\psi\in L^{\infty}(X\times X)$ be the
function given by
$\psi(x,y)=\sum_{k=1}^{\infty} a_k(x)b_k(y)$.
For all $f\in L^2(X\times X)$ we have
$S_\nph(T_f)=T_{\nph f}$ and $\sum_{k=1}^{\infty}M_{b_k}T_fM_{a_k}=T_{\psi f}$.
We obtain that $\nph(x,y)
f(x,y)=\psi(x,y) f(x,y)$ for every $f\in L^2(X\times X)$
and this implies that $\nph(x,y)=\psi(x,y)$ for almost all $(x,y)\in X\times X$.

To obtain the converse statement one should first note that for
$\nph(x,y) = \sum_{k=1}^{\infty} a_k(x)b_k(y)$ we have that the operator
$S_{\nph}(T)$ is given by (\ref{sum}), for all $T\in \cl C_2(L^2(X))$.
To complete the proof, one can apply estimation arguments similar to the one in (\ref{estim}).
\prend

Certain notions pertinent to Measure Theory that were initially introduced
by Arveson \cite{a} and later developed in \cite{eks} play
an important role in the study of measurable Schur multipliers.
A subset $E\subseteq X\times Y$ is called {\bf marginally null} if $E\subseteq (M\times Y) \cup (X\times N)$
for some null sets $M\subseteq X$ and $N\subseteq Y$. Two measurable sets $E,F\subseteq X\times Y$ are called
{\bf marginally equivalent} if their symmetric difference is marginally null.
A measurable subset $\kappa\subseteq X\times Y$ is called {\bf $\omega$-open} if $\kappa$
is marginally equivalent to a subset of the form $\cup_{i=1}^{\infty}\alpha_i\times\beta_i$,
where $\alpha_i\subseteq X$ and $\beta_i\subseteq Y$ are measurable;
$\kappa$ is called {\bf $\omega$-closed} if
its complement is $\omega$-open.
The collection of all $\omega$-open subsets of $X\times Y$ is an {\bf $\omega$-topology};
that is, a family of sets closed under finite intersections and countable unions, containing the
empty set and the set $X\times Y$. The described $\omega$-topology plays an important role
in the theory of masa-bimodules (see \cite{eks}). The morphisms of $\omega$-topological spaces are
$\omega$-continuous mappings. A function $f : X\times Y\rightarrow \bb{C}$ is {\bf $\omega$-continuous}
if $f^{-1}(U)$ is an $\omega$-open set for every open subset $U\subseteq \bb{C}$.
The following was shown in \cite{ks}:

\begin{proposition}
Every $\mu,\nu$-multiplier coincides almost everywhere with an
$\omega$-continuous function.
\end{proposition}

From this result, we immediately obtain that if $\chi_{\kappa}$ is
an idempotent Schur multiplier, 
where $\kappa\subseteq X\times Y$ then $\kappa$ is differs by a null
set from a subset that is both $\omega$-closed and $\omega$-open.
This fact was pointed out in \cite{kp}.

The predual of $\cl B(L^2(X),L^2(Y))$ can be naturally identified with the projective
tensor product $L^2(X)\hat{\otimes}L^2(Y)$.
Suppose that $h\in L^{\infty}(X)\hat{\otimes}L^2(Y)$ and let $h =
\sum_{k=1}^{\infty} f_k\otimes g_k$ be an associated series for $h$,
where $\sum_{k=1}^{\infty} \|f_k\|_2^2 < \infty$ and
$\sum_{k=1}^{\infty} \|g_k\|_2^2 < \infty$. These conditions easily
imply that the formula
$$h(x,y) = \sum_{k=1}^{\infty} f_k(x)g_k(y), \ \ \ (x,y)\in X\times Y,$$
defines a function (which we denote again by $h$). If the $f_k$'s or the $g_k$'s
are replaced by some equivalent functions (with respect to the measures $\mu$ and $\nu$)
then the function $h$ will change only on a marginally null set.

Let $\Gamma(X,Y)$ be the set of all functions $h$ defined as above.
Another useful result of Peller \cite{peller_two_dim} identifies
$S(X,Y)$ as the multiplier algebra of $\Gamma(X,Y)$:

\begin{proposition}
A function $\nph\in L^{\infty}(X\times Y)$ is a Schur multiplier if and only if
$\nph h \in \Gamma(X,Y)$ for every $h\in \Gamma(X,Y)$.
\end{proposition}
\proof For $k\in L^2(X\times Y)$ and $h\in \Gamma(X,Y)$, one has
$\langle T_k,T_h\rangle = \int kh d(\mu\times \nu)$. It follows that
$h\in \Gamma(X,Y)$ if and only if $$\exists C > 0 \mbox{ with } \left|\int kh d(\mu\times \nu)\right|\le
 C \|T_k\|_{\op} \text{ for all }k\in L^2(X\times Y).$$ Now we have that if
$\nph$ is a Schur multiplier then
\begin{eqnarray*}
\left|\int \varphi hk d(\mu\times \nu)\right| = \left|\langle
S_{\varphi}(T_k), T_h \rangle \right| \le
|\nph\|_{\mm}\|T_h\|_1\|T_k\|_{\op}  \Leftrightarrow \varphi h\in
\Gamma(X,Y).
\end{eqnarray*}
Conversely, if $\nph h\in \Gamma(X,Y)$  for all $h\in\Gamma(X,Y)$,
then by the Closed Graph Theorem, the mapping $T_h\mapsto T_{\nph
h}$ on $\cl C_1(L^2(X), L^2(Y))$ is bounded and
$$|\langle S_{\varphi}(T_k), T_h \rangle | =|\int \varphi hk
d(\mu\times \nu)|=|\langle T_k, T_{\nph h} \rangle |\leq
C\|T_k\|_{\op}\|T_{h}\|_1.$$ Thus, $\|S_{\varphi}(T_k)\|_{\op}\leq C\|T_k\|_{\op}$.
\prend

Measurable Schur multipliers are closely related to the theory of double operator integrals
developed by Birman and Solomyak in a series of papers \cite{BS1}, \cite{BS2}, \cite{BS3}
(see also the survey paper \cite{BS4}). We describe here briefly this connection.
Let $E(\cdot)$ and $F(\cdot)$ be spectral measures defined on measure spaces
$X$ and $Y$, and taking values in the projection lattices
of Hilbert spaces $H$ and $K$, respectively.
We fix a scalar valued measure $\mu$ (resp. $\nu$) on $X$ (resp. $Y$), equivalent to $E$ (resp. $F$).
A {\bf double operator integral}
is a formal expression of the form
\begin{equation}\label{doi}
\cl I^{E,F}_{\nph}(T) = \int \nph(x,y) d E(x) T d F(y),
\end{equation}
where $\nph$ is an essentially bounded function defined on $X\times Y$,
and $T\in \cl B(K,H)$. A precise meaning can be given to (\ref{doi}) as follows.
Let $G$ be the (unique) spectral measure defined on the product $X\times Y$ of the
measure spaces $X$ and $Y$ and taking values in the projection lattice of the
Hilbert space $\cl C_2(K,H)$, given on measurable rectangles by
$G(\alpha\times \beta)(T) = E(\alpha)TF(\beta)$.
For $\nph\in L^{\infty}(X\times Y, \mu\times \nu)$ and $T\in \cl C_2(K,H)$, let
$$\cl I^{E,F}_{\nph}(T)\stackrel{def}{=} \left(\int_{X\times Y} \nph(x,y) d G(x,y)\right)(T).$$
Spectral Theory gives a precise meaning to the right hand side of
the above expression, a well-defined element of $\cl C_2(K,H)$. Note
that $\cl I^{E,F}_{\nph}(T)$ depends only on the equivalence class
(with respect to product measure) of the function $\nph$. For some
functions $\nph$, there may exist a constant $C > 0$ such that
\begin{equation}\label{doi2}
\|\cl I^{E,F}_{\nph}(T)\|_{\op}\leq C\|T\|_{\op};
\end{equation}
in this case the mapping
$\cl I^{E,F}_{\nph}$ extends uniquely to a bounded
mapping on $\cl K(K,H)$ and, after taking second duals,
to a bounded mapping on $\cl B(K,H)$.

Suppose for a moment that the spectral measures $E$ and $F$ are multiplicity free, that is, the
(abelian) von Neumann algebras generated by their ranges are maximal. It is easy to note that in this case
a function $\nph$ satisfies the inequality (\ref{doi2}) for some constant $C > 0$ if and only
if $\nph$ is a Schur multiplier.

If $E$ possesses non-trivial multiplicity then one can decompose the Hilbert space $H$
as a direct integral $H = \int_X H(x) d\mu(x)$ in such a way that the elements of the abelian
von Neumann algebra generated by $E(\cdot)$ are precisely the diagonal operators. Similarly,
$K = \int_Y K(y) d\nu(y)$. The space $\cl C_2(K,H)$ possesses a decomposition as
$\cl C_2(K,H) = \int_{X\times Y} \cl C_2(K(y),H(x)) d\mu\times \nu(x,y)$.
Thus, every element $T$ of $\cl C_2(K,H)$ gives rise to an operator valued \lq\lq kernel''
$(T(x,y))_{x,y}$. We have that $\cl I^{E,F}_{\nph}(T)$ is the operator with corresponding
kernel $(\nph(x,y)T(x,y))_{x,y}$. In this way, double operator integrals may again be realised as
measurable Schur multipliers.

\bigskip

One application of double operator integrals is to perturbation
theory, in particular, to the study of operators of the form
$h(A)-h(B)$ depending on the properties of $A-B$, where $A$, $B$ are
selfadjoint operators and $h$ is a function defined on an interval
containing the spectra of $A$ and $B$.

Assume $h(x)$ is a uniformly Lipschitz function on $\mathbb R$, and
let
$$\check h(x,y)=\frac{h(x)-h(y)}{x-y}.$$
This function is well-defined and bounded outside the diagonal
$\{(x,x) : x\in \bb{R}\}$. Since the diagonal has measure zero,
$\check h(x,y)$ is well defined up to a null set. Thus, extending it
in an arbitrary fashion to a bounded function defined on the whole
of $\bb{R}\times\bb{R}$ always yields the same element of
$L^{\infty}(\bb{R}\times\bb{R})$ which we again denote by
$\check{h}$.

\begin{theorem}\cite{BS4} Let $E(\cdot)$, $F(\cdot)$ be the spectral measures of $A$ and $B$ respectively.
Suppose that $\check h$ satisfies (\ref{doi2}). Then 
$$h(A)-h(B)=\cl I^{E,F}_{\check h}(A-B)$$
and hence
$$\|h(A)-h(B)\|\leq C\|A-B\|.$$
\end{theorem}

The property of $\check h$ being a Schur multiplier is closely
related to a kind of "operator smoothness" of $h$.
Let $\alpha$ be a compact set in ${\mathbb R}$. A continuous
function $h$ on $\alpha$ is called Operator Lipschitz  on $\alpha$
if there is $D>0$ such that
$$\|h(A)-h(B)\|\leq D\|A-B\|$$
for all selfadjoint operators $A$, $B$ with spectra in $\alpha$.

\begin{proposition}\label{eqdiv}
Let $I$ be a compact interval. A function $\check h$
is a Schur multiplier on $I\times I$ if and only if $h$ is Operator
Lipschitz on $I$.
\end{proposition}

The proof of this result is based on the result by Kissin and Shulman
\cite{ks_edinburgh} that for all compact sets $I$ of $\mathbb R$ a
function $h$ is Operator Lipschitz on $I$ if and only if $h$ is
commutator bounded, that is, there exists $D>0$ such that for any selfadjoint
$A$ with spectrum in $I$ and any bounded operator $X$, the inequality
\begin{equation}\label{eq_eq}
\|h(A)X-Xh(A)\|\leq D\|AX-XA\|
\end{equation}
holds. Note also that if
$A$ is the multiplication operator $(Ag)(x)=xg(x)$ on $L^2(I,\mu)$
and $X = T_k$, where $k(x,y)=(x-y)k_1(x,y)$ for some
$k_1\in L^2(I\times I,\mu\times\mu)$, then (\ref{eq_eq}) is equivalent
to
$$\|S_{\check h}(T_k)\|\leq D\|T_k\|.$$

For other applications, in particular to differentiation of
functions of selfadjoint operators, we refer the reader to
\cite{BS4}, and for more results on applications of Schur
multipliers to operator inequalities -- to \cite{HK}.

\bigskip

As the connections with double operator integrals shows, the purpose
of introducing measurable Schur multipliers is not limited to
studying the notion in the greatest possible generality. We now
further illustrate this, describing a connection with Harmonic
Analysis.
Let $G$ be a locally compact group which we assume for technical simplicity to
be $\sigma$-compact. We let $L^p(G)$, $p = 1,2,\infty$, be the
corresponding function spaces with respect to the (left) Haar measure.
By $\lambda : G\rightarrow \cl B(L^2(G))$ we denote the left regular representation
of $G$; thus, $\lambda_s f(t) = f(s^{-1}t)$, $s,t\in G$, $f\in L^2(G)$.
We recall that the Fourier algebra $A(G)$ of $G$ is the space of
all \lq\lq matrix coefficients of $G$ in its left regular representation'', that is,
$$A(G) = \{s\rightarrow (\lambda_s\xi,\eta) : \xi,\eta\in L^2(G)\}.$$
If $G$ is commutative, then $A(G)$ is the image of
$L^1(\hat G)$ under Fourier transform, where $\hat G$ is the dual group of $G$. The Fourier
algebra of general locally compact groups was introduced and
studied (along with other objects pertinent to Non-commutative
Harmonic Analysis) by Eymard in \cite{eymard}. It is a commutative
regular semi-simple Banach algebra of continuous functions vanishing
at infinity and has $G$ as its spectrum. Moreover, its Banach space dual is
isometric to the von Neumann algebra $\vn(G)$ of $G$, that is, the
weakly closed subalgebra of $\cl B(L^2(G))$ generated by the
operator $\lambda_s$, $s\in G$. The duality between these two spaces
is given by the formula $\langle \lambda_x,f\rangle = f(x)$ (recall that
$f\in A(G)$ is a function on $G$).

A function $\nph\in L^{\infty}(G)$ is called a multiplier of $A(G)$ if
$\nph f\in A(G)$ for every $f\in A(G)$. The identification of the multipliers of $A(G)$
(which, as is easy to see, form a function algebra on $G$)
has received a considerable attention in the literature. A classical result in this direction
(\cite{rudin}) states that if $G$ is abelian then $\nph$ is a multiplier of $A(G)$
if and only if it is the Fourier transform of a regular Borel measure on $\hat{G}$.

Up to date, a satisfactory characterisation of the multiplier algebra of $A(G)$ is not known.
There is, however, a neat and useful characterisation of the
subalgebra $M_{cb}A(G)$ of {\bf completely bounded multipliers of $A(G)$}
introduced in \cite{ch}.
In order to define these multipliers, we need the notions of an operator space and
a completely bounded map.

An {\bf operator space} is a complex vector space $\cl X$ for which
a norm $\|\cdot\|_n$ is given on the space $M_n(\cl X)$ of all $n$
by $n$ matrices with entries in $\cl X$ satisfying the following
conditions, called Ruan's axioms:

\smallskip

(R1) $\|x \oplus y\|_{n+m} = \max\{\|x\|_n,\|y\|_m\}$, for all
$m,n\in \bb{N}$ and all $x\in M_n(\cl X)$ and $y\in M_m(\cl X)$;

\smallskip

(R2) $\|\alpha \cdot x \cdot \beta\|_{m} \leq \|\alpha\|\|\beta\|\|x\|_n$, for all $x\in M_n(\cl X)$,
$\alpha\in M_{m,n}(\bb{C})$, $\beta \in M_{n,m}(\bb{C})$ and all $m,n\in \bb{N}$.

\smallskip

\noindent In property (R2), we have denoted by $\cdot$ the natural
left action of $M_{m,n}$ on $M_n(\cl X)$ as well as the natural
right action of $M_{n,m}$ on $M_n(\cl X)$. By $\|\alpha\|$ we mean
the norm of the matrix $\alpha\in M_{m,n}$ when considered as an
operator from $\bb{C}^n$ into $\bb{C}^m$ (where $\bb{C}^n$ and
$\bb{C}^m$ are equipped with the $\ell^2$-norm).

If $\cl X$ and $\cl Y$ are operator spaces and  $\phi : \cl
X\rightarrow\cl Y$ is a linear map, then one may consider the maps
$\phi_n : M_n(\cl X)\rightarrow M_n(\cl Y)$ defined by applying the
mapping $\phi$ entry-wise. The map $\phi$ is called {\bf completely
bounded} if each $\phi_n$, $n\in \bb{N}$, is bounded, and
$\sup_{n\in \bb{N}} \|\phi_n\| < \infty$. It is called a {\bf
complete isometry} if $\phi_n$ is an isometry for every $n\in
\bb{N}$. Since complete isometries preserve all the given structure
of an operator space, they constitute the right notion of an
isomorphism in the category of operator spaces.

If $\cl X$ is an operator space and $\cl X^*$ is its dual Banach
space, then the space $M_n(\cl X^*)$ can be naturally identified
with the space $CB(\cl X, M_n(\bb{C}))$ of all completely bounded
linear maps from $\cl X$ into $M_n(\bb{C})$. If we equip $M_n(\cl
X^*)$ with the norm coming from this identification, then the
obtained family of norms turns $\cl X^*$ into an operator space.

If $\cl X\subseteq \cl B(H)$ for some Hilbert space $H$, then
$M_n(\cl X)$ is naturally embedded into $\cl B(H^n)$ and hence inherits its
operator norm. The family of norms obtained in this way satisfy Ruan's axioms.
Ruan's representation theorem asserts that every operator space is completely isometric
to a subspace of $\cl B(H)$, for some Hilbert space $H$.

After this very brief introduction of the basic notions of Operator
Space Theory, we can return to multipliers. From the last paragraph
it follows that every von Neumann algebra, in particular $\vn(G)$,
is an operator space. Hence its dual, and therefore its predual
$A(G)$ (which is a subspace of its dual) possesses a canonical
operator space structure. A multiplier $\nph$ of $A(G)$ is called
completely bounded if the mapping $f\rightarrow \nph f$ is
completely bounded. Namely, the following was established by Bozejko
and Fendler in \cite{bf} (see also \cite{Pi}):

\begin{theorem}\label{bf}
If $\nph\in L^{\infty}(G)$, let $\tilde{\nph}\in L^{\infty}(G\times G)$ be the
function given by $\tilde{\nph}(s,t) = \nph(s^{-1}t)$.
A function $\nph\in L^{\infty}(G)$ belongs to $M_{cb}A(G)$ if and only if
the function $\tilde{\nph}$ is a Schur multiplier with respect to the left Haar measure.
\end{theorem}

\section{Going non-commutative}\label{sec3}

Let $(X,\mu)$ and $(Y,\nu)$ be standard ($\sigma$-finite) measure
spaces. It is immediate from the definitions that if $\nph,\psi\in
S(X,Y)$ then $S_{\nph}S_{\psi} = S_{\nph\psi}$. Thus,
$S_{\nph}S_{\psi} = S_{\psi}S_{\nph}$ for every $\nph,\psi\in
S(X,Y)$; in other words, the collection of all mappings $\{S_{\nph}
: \nph\in S(X,Y)\}$ is commutative. In view of the contemporary
trend in Functional Analysis to seek non-commutative versions of
\lq\lq classical'' notions and results \cite{blm}, \cite{er},
\cite{paulsen}, \cite{pisier_intr}, it is natural to ask whether
there is a natural non-commutative, or \lq\lq quantised'' version of
Schur multipliers. This question was pursued by Kissin and Shulman
in \cite{ks}, and is the topic of this section.

For a Hilbert space $H$, we write $H^{\dd}$ for the (Banach space) dual of $H$. There
exists a conjugate-linear surjective isometry $\partial : H\rightarrow H^{\dd}$
given by $\partial(x) (y) = (y,x)$, $x,y\in H$.

Let $H$ and $K$ be Hilbert spaces and $\theta : H\otimes K
\rightarrow \cl C_2(H^{\dd},K)$ be the natural surjective isometry
from the Hilbert space tensor product $H\otimes K$ onto the space
$\cl C_2(H^{\dd},K)$ of all Hilbert-Schmidt operators from $H^{\dd}$
into $K$ given by $\theta(x\otimes y)(z^{\dd}) = (x,z)y$. This
identification allows us to equip $H\otimes K$ with an \lq\lq
operator'' norm: if $\xi\in H\otimes K$, let $\|\xi\|_{\op} =
\|\theta(\xi)\|_{\op}$. We call an element $\nph\in \cl B(H\otimes
K)$ a {\bf concrete operator multiplier} if there exists a constant
$C > 0$ such that $\|\nph\xi\|_{\op} \leq C \|\xi\|_{\op}$, for
every $\xi \in H\otimes K$. We call the smallest possible constant
$C$ with this property the concrete multiplier norm of $\nph$. It
follows from this definition that the set $\frak{M}(H,K)$ of all
concrete operator multipliers on $H\otimes K$ is a subalgebra of
$\cl B(H\otimes K)$. It is also immediate that if $H = L^2(X,\mu)$
and $K = L^2(Y,\nu)$ for some standard measure spaces $(X,\mu)$ and
$(Y,\nu)$ and $\nph$ is the multiplication operator on $L^2(X\times
Y) = L^2(X)\otimes L^2(Y)$ corresponding to a function
$\tilde{\nph}\in L^{\infty}(X\times Y)$ then $\nph$ is a concrete
operator multiplier if and only if $\tilde{\nph}$ is a Schur
$\mu,\nu$-multiplier. Thus, the algebra $\frak{M}(H,K)$ contains
$S(X,Y)$ as a commutative subalgebra. Note that there are \lq\lq
many'' commutative subalgebras of $\frak{M}(H,K)$ of this form, one
for each realisation
of $H$ and $K$ as $L^2$-spaces over some standard measure spaces.

In the commutative theory, a special interest has been paid to the case where
the measure spaces $(X,\mu)$ and $(Y,\nu)$ are regular Borel spaces of
complete metrisable topologies, and the multipliers $\nph$ are
continuous functions on $X\times Y$. The non-commutative expression of
continuous functions is given in terms of C*-algebras. Therefore, it is natural to
extend the setting of concrete operator multipliers given above as follows.
Suppose that $\cl A$ and $\cl B$ are unital C*-algebras and
$\pi : \cl A\rightarrow \cl B(H_{\pi})$ and $\rho : \cl B\rightarrow \cl B(H_{\rho})$
are *-representations.
It is well-known that there exists a unique *-representation $\pi\otimes\rho$ of the
minimal tensor product $\cl A\otimes\cl B$ of $\cl A$ and $\cl B$ on
$H_{\pi}\otimes H_{\rho}$.
An element $\nph\in \cl A\otimes \cl B$ is called a
{\bf $\pi,\rho$-multiplier} \cite{ks} if
$(\pi\otimes\rho)(\nph)$ is a concrete operator multiplier.
We let $\|\nph\|_{\pi,\rho}$ be the concrete operator multiplier norm of $(\pi\otimes\rho)(\nph)$.

Let $\frak{M}^{\pi,\rho}(\cl A,\cl B)$ be the set of all
$\pi,\rho$-multipliers. It is immediate that
$\frak{M}^{\pi,\rho}(\cl A,\cl B)$ is a subalgebra of $\cl A\otimes\cl B$
containing the algebraic tensor product $\cl A\odot \cl B$.
To see the last statement, note that if $a\in \cl A$ and $b\in \cl B$
then
$\theta((\pi\otimes\rho) (a\otimes b)(\xi)) = \rho(b)\theta(\xi)\pi(a)^{\dd}$,
for every $\xi\in H_{\pi}\otimes H_{\rho}$; we hence have that
$\|a\otimes b\|_{\pi,\rho}\leq \|a\|\|b\|$.

We let $\frak{M}(\cl A,\cl B) = \cap_{\pi,\rho} \frak{M}^{\pi,\rho}(\cl A,\cl B)$
where the intersection is taken over all representations $\pi$ of $\cl A$ and $\rho$
of $\cl B$.
The elements of $\frak{M}(\cl A,\cl B)$ are called {\bf universal multipliers}. By the previous paragraph,
every element of the algebraic tensor product $\cl A\odot\cl B$ is a universal multiplier.
It is not difficult to see that if $\nph\in \frak{M}(\cl A,\cl B)$ then
$\|\nph\|_{\mm} \stackrel{def}{=} \sup_{\pi,\rho}\|\nph\|_{\pi,\rho}$ is finite; we call
$\|\nph\|_{\mm}$ the (universal) multiplier norm of $\nph$.

Two immediate questions arise:

\medskip

(a) How does the algebra $\frak{M}^{\pi,\rho}(\cl A,\cl B)$
depend on $\pi$ and $\rho$?

\medskip

(b) Is there a characterisation of its elements extending the
Grothendieck-Peller's characterisation of Schur multipliers?

\medskip

It was observed by Kissin
and Shulman \cite{ks} that
 (a) is related to the notion of approximate equivalence of
representations due to Voiculescu \cite{voiculescu} and its
extension, the approximate sub-ordinance introduced by Hadwin
\cite{hadwin}.  We recall these notion here.
Let $\pi$ and $\pi'$ be
$*$-representations of a $C^*$-algebra ${\cl A}$ on Hilbert spaces
$H$ and $H'$, respectively. We say that $\pi'$ is {\it approximately
subordinate} to $\pi$ and write $\pi'\stackrel{a}\ll\pi$ if there is
a net $\{U_{\lambda}\}$ of isometries from $H'$ to $H$ such that
\begin{equation}\label{as}
\|\pi(a)U_{\lambda}-U_{\lambda}\pi'(a)\|\to 0\text{ for all
}a\in{\cl A}.
\end{equation}
The representations $\pi'$ and $\pi$ are said to be {\it
approximately equivalent} if the operators $U_{\lambda}$ can be
chosen to be unitary; in this case we write
$\pi'\stackrel{a}\sim\pi$.

The following result was established in \cite{ks}:

\begin{theorem}[Comparison Theorem]\label{compar}
Let $\cl A$ and $\cl B$ be C*-algebras and $\pi$, $\pi'$ (resp. $\rho$, $\rho'$) be
representations of $\cl A$ (resp. $\cl B$). Suppose that $\pi'\stackrel{a}\ll\pi$ and
$\rho'\stackrel{a}\ll\rho$. Then
$\frak{M}^{\pi,\rho}(\cl A,\cl B)\subseteq \frak{M}^{\pi',\rho'}(\cl A,\cl B)$.
Moreover, if $\nph\in \frak{M}^{\pi,\rho}(\cl A,\cl B)$ then
$\|\nph\|_{\pi',\rho'}\leq \|\nph\|_{\pi,\rho}$.

In particular, if $\pi'\stackrel{a}\sim\pi$ and $\rho'\stackrel{a}\sim\rho$
then
$\frak{M}^{\pi,\rho}(\cl A,\cl B)=\frak{M}^{\pi',\rho'}(\cl A,\cl B)$ and
$\|\nph\|_{\pi',\rho'} = \|\nph\|_{\pi,\rho}$ for every $\nph\in \frak{M}^{\pi,\rho}(\cl A,\cl B)$.
\end{theorem}

We note that, by \cite{hadwin}, $\pi'\stackrel{a}\ll\pi$ if and only if
$$\rank\pi'(a) \leq \rank \pi(a), \ \ \ \mbox{for every } a\in \cl A.$$

Theorem \ref{compar} has some interesting consequences about measurable and classical Schur multipliers \cite{ks}:

\begin{corollary}\label{cor_supp}
Let $X$ and $Y$ be locally compact Hausdorff spaces with countable
bases, and let $\mu$ and $\mu'$ (resp. $\nu$ and $\nu'$) be
$\sigma$-finite Borel measures on $X$ (resp. $Y$). Suppose that
$\supp \mu'\subseteq \supp \mu$ and $\supp \nu'\subseteq \supp\nu$.
Then every $\mu,\nu$-multiplier in $C_0(X\times Y)$ is also a
$\mu',\nu'$-multiplier.

In particular, if $\supp \mu = X$ and $\supp \nu = Y$ then an element
$\nph\in C_0(X\times Y)$ is a $\mu,\nu$-multiplier if and only if
$\nph$ is a classical Schur multiplier on $X\times Y$.
\end{corollary}

There is a version of the last result for functions $\nph$ that are not required to be continuous \cite{ks}:

\begin{theorem}
Let $(X,\mu)$ and $(Y,\nu)$ be standard $\sigma$-finite measure spaces and $\nph\in L^{\infty}(X\times Y)$
be an $\omega$-continuous function. The following are equivalent:

(i) \ $\nph$ is a $\mu,\nu$-multiplier;

(ii) there exist null sets $X_0\subseteq X$ and $Y_0\subseteq Y$ such that
the restriction $\tilde{\nph}$ of $\nph$ to $(X\setminus X_0)\times (Y\setminus Y_0)$ is a classical
Scur multiplier.

Moreover, if (i) holds then the sets $X_0$ and $Y_0$ can be chosen
in such a way that the $\mu,\nu$-multiplier norm of $\nph$ equals
the classical Schur multiplier norm of $\tilde{\nph}$.
\end{theorem}

We now address Question (b) above concerning the characterisation of
operator multipliers. At the moment, no such characterisation is
known for the classes $\frak{M}^{\pi,\rho}(\cl A,\cl B)$. The reason
lies in the lack of complete boundedness, which we now explain.
Suppose that $\cl A\subseteq \cl B(H)$ and $\cl B\subseteq \cl B(K)$
are concrete C*-algebras, and take an element $\nph\in
\frak{M}^{\id,\id}(\cl A,\cl B)$, where $\id$ denotes the identity
representations of $\cl A$ and $\cl B$ on $H$ and $K$, respectively.
Since $\nph$ is a concrete operator multiplier on $H\otimes K$, we
have that the mapping $S_{\nph} : \cl C_2(H^{\dd},K)\rightarrow \cl
C_2(H^{\dd},K)$ given by $S_{\nph}(\theta(\xi)) = \theta(\nph\xi)$,
has a canonical extension to a bounded mapping (denoted in the same
way) $S_{\nph} : \cl K(H^{\dd},K)\rightarrow \cl K(H^{\dd},K)$. By
passing to second duals, we arrive at a bounded mapping
$S_{\nph}^{**} : \cl B(H^{\dd},K)\rightarrow \cl B(H^{\dd},K)$. In
general, however, the mappings $S_{\nph}$ and $S_{\nph}^{**}$ need
not be completely bounded. If $\nph$ is assumed to lie in the
smaller class $\frak{M}(\cl A,\cl B)$ of universal multipliers, then
the mappings $S_{\nph}$ and $S_{\nph}^{**}$ turn out to be
completely bounded. This can be seen by considering the $n$-fold
ampliations $\id^{(n)}$ of the identity representations of $\cl A$
and $\cl B$. In fact, the following statement \cite{jtt} shows that
in order to decide whether an element $\nph\in \cl A\otimes\cl B$ is
a universal multiplier, it suffices to check that it is a
$\id^{(n)},\id^{(n)}$-multiplier for all $n\in \bb{N}$. The proof of
this result uses the Comparison Theorem \ref{compar}.

\begin{proposition}\label{cb}
Let $\cl A\subseteq \cl B(H)$ and $\cl B\subseteq \cl B(K)$.
An element $\nph\in \cl A\otimes\cl B$ is a universal multiplier if and only if the
(weak* continuous) mapping $S_{\nph}^{**}$ is completely bounded.
\end{proposition}

The structure of normal (that is, weak* continuous) completely
bounded maps on $\cl B(K,H)$ is well understood. We recall a
well-known result of Haagerup \cite{haagerup}: a mapping $\Phi : \cl
B(K,H)\rightarrow\cl B(K,H)$ is normal and completely bounded if and
only if there exist families $\{a_k\}_{k=1}^{\infty}\subseteq \cl
B(H)$ and $\{b_k\}_{k=1}^{\infty}\subseteq \cl B(K)$ such that the
series $\sum_{k=1}^{\infty} a_k a_k^*$ and $\sum_{k=1}^{\infty}
b_k^* b_k$ are weak* convergent and
$$\Phi(x) = \sum_{k=1}^{\infty} a_k x b_k, \ \ \ \ \mbox{for all } x\in \cl B(K,H).$$
Thus, one may associate with every normal completely bounded map $\Phi$ on
$\cl B(K,H)$ a formal series $\sum_{k=1}^{\infty} a_k\otimes b_k$ where the
families $\{a_k\}_{k=1}^{\infty}\subseteq \cl B(H)$ and
$\{b_k\}_{k=1}^{\infty}\subseteq \cl B(K)$ are assumed to satisfy the above convergence
conditions. Two such formal series are identified if the corresponding mappings are equal.
The collection of all such series is known as the weak* (or the extended) Haagerup tensor product of
$\cl B(H)$ and $\cl B(K)$ and denoted by $\cl B(H)\otimes_{\eh} \cl B(K)$. In fact,
$\cl B(H)\otimes_{\eh} \cl B(K)$ can be viewed as a certain weak completion of the algebraic
tensor product $\cl B(H)\odot \cl B(K)$ (see \cite{bs} where this tensor product was introduced).

The weak* Haagerup tensor product can be defined for any pair of
dual operator spaces $\cl X^*$ and $\cl Y^*$ and is the dual
operator space of the {\bf Haagerup tensor product} $\cl
X\otimes_{\hh}\cl Y$ of $\cl X$ and $\cl Y$. The latter is the
tensor product that linearises completely bounded bilinear mappings.
These are defined as follows: Suppose that $\phi : \cl X\times\cl
Y\rightarrow\cl Z$ is a bilinear mapping, where $\cl Z$ is another
operator space. One may define the mappings $\phi^{(n)} : M_n(\cl
X)\times M_n(\cl Y)\rightarrow M_n(\cl Z)$, $n\in \bb{N}$, by
$\phi^{(n)} ((x_{i,j}),(y_{k,l})) = (\sum_{k=1}^n
\phi(x_{i,k},y_{k,j})_{i,j}$. The mapping $\phi$ is called
completely bounded if $\|\phi\|_{cb} \stackrel{def}{=}\sup_{n\in
\bb{N}}\|\phi^{(n)}\| < \infty$. The operator space $\cl
X\otimes_{\hh}\cl Y$ has the property that for every completely
bounded bilinear mapping $\phi : \cl X\times\cl Y\rightarrow\cl Z$
the linearised mapping $\tilde{\phi}$ is completely bounded as a map
from $\cl X\otimes_{\hh}\cl Y$ into $\cl Z$ with the same cb norm.
We note that we will introduce later a generalisation of the above
notion of complete boundedness to multilinear maps.

The correspondence between elements of $\cl B(H)\otimes_{\eh}\cl B(K)$ and
normal completely bounded mappings on $\cl B(K,H)$ is bijective; if $u\in \cl B(H)\otimes_{\eh}\cl B(K)$,
we write $\Phi_u$ for the corresponding mapping.
The space $\cl B(H)\otimes_{\eh}\cl B(K)$ is an operator space in its own right
(indeed, we have the completely isometric identification $\cl B(H)\otimes_{\eh}\cl B(K)
= (\cl C_1(H)\otimes_{\hh}\cl C_1(K))^*$), and the norm $\|u\|_{\eh}$ of
an element $u$ is equal to the completely bounded norm of $\Phi_u$.

The extended Haagerup tensor product can be defined for every pair of operator spaces \cite{effros_ruan}.
To do this, we follow the approach given in \cite{spronk}.
Let $\cl E\subseteq \cl B(H)$ and $\cl F\subseteq \cl B(K)$ be {\it norm closed} subspaces.
Then the extended Haagerup tensor product $\cl E\otimes_{\eh}\cl F$ of $\cl E$ and $\cl F$ is
the subspace
$$\{u\in \cl B(H)\otimes_{\eh}\cl B(K) : \id\otimes\omega (u)\in \cl E, \tau\otimes\id (u)\in \cl F, \
\forall \ \omega\in \cl B(K)_*, \tau \in \cl B(H)_*\}.$$
Here, $\id\otimes\omega$ (resp. $\tau\otimes\id$)
is the left (resp. right) slice map from $\cl B(H)\otimes_{\eh}\cl B(K)$ into $\cl B(H)$
(resp. $\cl B(K)$) along the functional $\omega$ (resp. $\tau$). The fact that these maps are well-defined
needs a justification that we omit.

The extended Haagerup tensor product of operator spaces is functorial: if
$f : \cl E\rightarrow\cl E'$ and $g : \cl F\rightarrow\cl F'$ are completely bounded
maps then there exists a unique completely bounded map $f\otimes g : \cl E\otimes_{\eh}\cl F
\rightarrow \cl E'\otimes_{\eh}\cl F'$ given on elementary tensors by
$(f\otimes g)(a\otimes b) = f(a)\otimes g(b)$. Moreover, if $f$ and $g$ are complete isometries
then so is $f\otimes g$ \cite{effros_ruan}.

We now return to universal multipliers. Recall that we have fixed
two C*-algebras $\cl A\subseteq \cl B(H)$ and $\cl B\subseteq \cl
B(K)$. Suppose that $\nph\in\cl A\otimes\cl B$ is a universal
multiplier. By Proposition \ref{cb}, the map $S_{\nph}^{**} : \cl
B(H^{\dd},K)\rightarrow\cl B(H^{\dd},K)$ is completely bounded.
Thus, by Haagerup's result described above, there exists an element
$u\in \cl B(K)\otimes_{\eh}\cl B(H^{\dd})$ such that $S_{\nph}^{**}
= \Phi_u$. It was shown in \cite{jltt} that, in fact, $u$ lies in
the extended Haagerup tensor product $\cl B\otimes_{\eh}\cl
A^{\dd}$, and that it does not depend on the concrete
representations of the C*-algebras $\cl A$ and $\cl B$ that we
started with. More precisely, the following \lq\lq symbolic
calculus'' result holds. (We denote by $\cl A^o$ the {\it opposite}
C*-algebra of $\cl A$ which coincides with $\cl A$ as an involutive
normed linear space but is equipped with the product $a\circ b =
ba$. For a representation $\pi : \cl A\rightarrow \cl B(H)$ we let
$\pi^{\dd} : \cl A^o \rightarrow \cl B(H^{\dd})$ be the
representation given by $\pi^{\dd}(a^o) = \pi(a)^{\dd}$.)

\begin{theorem}[Symbolic Calculus for universal multipliers]
Let $\cl A$ and $\cl B$ be C*-algebras. There exists an injective homomorphism
$\nph\rightarrow u_{\nph}$ from $\frak{M}(\cl A,\cl B)$ into $\cl B\otimes_{\eh}\cl A^o$
with the following universal property:
if $\pi : \cl A\rightarrow\cl B(H)$ and $\rho : \cl B\rightarrow \cl B(K)$ are *-representations
then
$$S_{\pi\otimes\rho(\nph)}^{**} = \Phi_{\rho\otimes\pi^{\dd}(u_{\nph})}.$$
Moreover, $\|\nph\|_{\mm} = \|u_{\nph}\|_{\eh}$, and
$u_{a\otimes b} = b\otimes a^o$, for all $a\in \cl A$, $b\in \cl B$.
\end{theorem}

We call the element $u_{\nph}$ the {\bf symbol} of the universal
multiplier $\nph$. The term \lq\lq symbolic calculus'' was first
used in the context of Schur multipliers by Katavolos and Paulsen
\cite{kp} where they explored the correspondence between a
measurable Schur multiplier $\nph$ and the mapping $S_{\nph}$.

Suppose now that $\nph\in \frak{M}(\cl A,\cl B)$. The corresponding
symbol $u_{\nph}$ has an associated series
$\sum_{i=1}^{\infty}b_i\otimes a_i^o$, where $a_i\in \cl A$ and
$b_i\in \cl B$, $i\in \bb{N}$. Let $u_N = \sum_{i=1}^{N}b_i\otimes
a_i^o$ and $\nph_N = \sum_{i=1}^{N}a_i\otimes b_i\in \cl A\odot\cl
B$, $N\in \bb{N}$. By Symbolic Calculus, $u_{\nph_N} = u_N$;
moreover, $u_{\nph} = {\rm w}^*$-$\lim_{N\rightarrow\infty} u_N$. We
also have that $\|\nph_N\|_{\mm} = \|u_N\|_{\eh} \leq
\|u_{\nph}\|_{eh}$ for all $N\in \bb{N}$. These observations can be
used to prove the next result which is the appropriate
generalisation of Grothendieck's and Peller's theorems. Before its
formulation, we note that it is easy to see that, for an element $v
\in \cl B\odot\cl A^o$, we have
$$\|v\|_{\eh} = \inf\left\{\left\|\sum_{i=1}^k d_id_i^*\right\|^{\frac{1}{2}}\left\|\sum_{i=1}^k c_ic_i^*\right\|^{\frac{1}{2}} : v = \sum_{i=1}^k d_i\otimes c_i^o\right\}.$$
For en element $\psi = \sum_{i=1}^k c_i\otimes d_i\in \cl A\odot\cl B$, we let
$\|\psi\|_{\ph} = \|\sum_{i=1}^k d_i\otimes c_i^o\|_{\eh}$.

\begin{theorem}[Characterisation Theorem]\label{charac}
Let $\cl A$ and $\cl B$ be C*-algebras and $\nph\in \cl A\otimes\cl B$. The following statements are
equivalent:

(i) \ $\nph\in \frak{M}(\cl A,\cl B)$ and $\|\nph\|_{\mm} < C$;

(ii) There exists a net $\{\nph_{\alpha}\}\subseteq \cl A\odot\cl B$
such that $\|\nph_{\alpha}\|_{\ph} < C$ for all $\alpha$ and
$(\pi\otimes\rho)(\nph_{\alpha})\rightarrow_{\alpha}
(\pi\otimes\rho)(\nph)$ weakly, for every pair $\pi,\rho$ of
irreducible representations of $\cl A$ and $\cl B$.
\end{theorem}

In the case $\cl A = C(X)$ and $\cl B = C(Y)$ are (unital) C*-algebras
($X$ and $Y$ being compact Hausdorff spaces), we obtain the following fact as
a consequence of the Characterisation Theorem, which shows that it
does extend Theorem \ref{groth2} to the non-commutative case:
If $\mu$ and $\nu$ are regular Borel measures on $X$ and $Y$, respectively,
then a function $\nph\in C(X\times Y)$ is a $\mu,\nu$-multiplier if and only if
there exist families $\{a_i\}_{i=1}^{\infty}\subseteq C(X)$ and
$\{b_i\}_{i=1}^{\infty}\subseteq C(Y)$ such that, if $\nph_k\in C(X\times Y)$
is given by $\nph_k(x,y) = \sum_{i=1}^k a_i(x)b_i(y)$, $(x,y)\in X\times Y$,
then $\sup_{k\in \bb{N}}\|\nph_k\|_{\mm} < \infty$ and $\nph_k\rightarrow\nph$ pointwise
$\mu\times\nu$-almost everywhere.

\bigskip

We now turn our attention to a subclass of universal multipliers; in order to define
it we recall a notion of compactness of completely bounded maps
introduced by Saar in \cite{saar}.
Let $\cl E$ and $\cl F$ be operator spaces and $\Phi : \cl E\rightarrow\cl F$ be a
completely bounded map. One calls $\Phi$ {\bf completely compact}
if for every $\epsilon > 0$ there exists a finite dimensional subspace $\cl F_0\subseteq \cl F$
such that $\dist(\Phi^{(n)}(x),M_n(\cl F_0))  < \epsilon$, for every $x$ in the unit ball
of $M_n(\cl X)$, and for every $n\in \bb{N}$. Clearly, every completely compact map is compact.

The following characterisation of completely compact maps on $\cl
K(H)$ was obtained in \cite{saar} and later in \cite{jltt} using a
different method. 

\begin{theorem}\label{saarcc}
A completely bounded map $\Phi : \cl K(H)\rightarrow\cl K(H)$ is completely compact
if and only if there exist sequences $\{a_i\},\{b_i\}\subseteq \cl K(H)$ such that
the series $\sum_{i=1}^{\infty}b_ib_i^*$ and $\sum_{i=1}^{\infty}a_i^*a_i$ are
norm convergent and
$$\Phi(x) = \sum_{i=1}^{\infty} b_i x a_i, \ \ \ \ \mbox{ for all } x\in \cl K(H).$$
\end{theorem}

We call an element $\nph\in \frak{M}(\cl A,\cl B)$ a {\bf compact} (resp.
{\bf completely compact}) multiplier if there exist faithful representations
$\pi$ and $\rho$ of $\cl A$ and $\cl B$, respectively, such that
the mapping $S_{\pi\otimes\rho(\nph)}^{**}$ is compact (resp. completely compact).
It is clear that every completely compact multiplier is compact.
It was shown in \cite{jltt} that an element $\nph\in \frak{M}(\cl A,\cl B)$
is (completely) compact if and only if
the mapping $S_{\pi_a,\rho_a}^{**}$ is (completely) compact, where
$\pi_a$ (resp. $\rho_a$) is the reduced atomic representation of $\cl A$
(resp. $\cl B$). This is rather natural to expect: let us recall a result
of Ylinen concerning the compact elements of a C*-algebra.
An element $a\in \cl A$ is called {\bf compact} if the mapping
$x\rightarrow axa$ on $\cl A$ is compact. Ylinen \cite{ylinen} showed that
an element $a\in \cl A$ is compact if and only if there exists a faithful representation
$\pi$ of $\cl A$ such that $\pi(a)$ is a compact operator; moreover, this happens
if and only if $\pi_a(a)$ is a compact operator.

Let us denote by $\cl K(\cl A)$ the set of all compact elements of
$\cl A$; it is well-known that $\cl K(\cl A)$ is a closed two sided
ideal of $\cl A$. By virtue of Ylinen's result, $\cl K(\cl A)$ is
*-isomorphic to a C*-algebra of compact operators, and is hence
*-isomorphic to a $c_0$-direct sum of the form $\oplus^{c_0}_{j\in
J} \cl K(H_j)$, for some index set $J$ and some Hilbert spaces
$H_j$, $j\in J$.

In the theorem that follows, we view the Haagerup tensor product
$\cl E\otimes_{\hh}\cl F$ of two operator spaces $\cl E$ and $\cl F$
as sitting completely isometrically in their extended Haagerup tensor product
$\cl E\otimes_{\eh}\cl F$.

\begin{theorem}\label{th_ccm}
Let $\cl A$ and $\cl B$ be C*-algebras and $\nph\in \frak{M}(\cl A,\cl B)$.
The following statements are equivalent:

(i) \ \ $\nph$ is a completely compact multiplier;

(ii) \ $u_{\nph}\in \cl K(\cl B)\otimes_{\hh}\cl K(\cl A^o)$;

(iii) there exists a sequence $\{\nph_k\}_{k=1}^{\infty}\subseteq \cl K(\cl A)\odot\cl K(\cl B)$
such that $\|\nph - \nph_k\|_{\mm}\rightarrow_{k\rightarrow\infty} 0$.
\end{theorem}

Of course, it is natural to ask what happens if $\nph$ is only assumed to be a compact
multiplier. One may show \cite[Proposition 7.1]{jltt} that if $\nph\in \frak{M}(\cl A,\cl B)$
is a compact multiplier then $u_{\nph}\in \cl K(\cl B)\otimes_{\eh}\cl K(\cl A^o)$; however,
the converse statement fails.
We do not have at present a complete characterisation of the compact universal multipliers;
however, the following \lq\lq automatic complete compactness'' result holds:

\begin{theorem}\label{th_automcc}
Let $\cl A$ and $\cl B$ be C*-algebras. Assume that
$\cl K(\cl A) \simeq \oplus^{c_0}_{i\in J_1} M_{n_i}$ and
$\cl K(\cl B) \simeq \oplus^{c_0}_{j\in J_2} M_{m_j}$, where
$\sup_{i\in J_1} n_i$ and $\sup_{j\in J_2} m_j$ are both finite.
Then every compact multiplier $\nph\in \frak{M}(\cl A,\cl B)$ is
automatically completely compact.
\end{theorem}

We now see that Theorem \ref{th_ccm} generalises Hladnik's
description of compact Schur multipliers (Theorem \ref{hladnik})
since in this case, by Theorem \ref{th_automcc}, every compact
multiplier is automatically completely compact.

In the case that $\sup_{i\in J_1} n_i$ and $\sup_{j\in J_2} m_j$ are both infinite,
we were able to exhibit in \cite{jltt} a compact multiplier $\nph\in \frak{M}(\cl A,\cl B)$ that is not
completely compact. The construction is based on an example of Saar \cite{saar} of a
compact completely bounded map on $\cl K(H)$ which is not completely compact.

\bigskip

In the last result that we mention in this section,
an answer is provided of when every universal multiplier
is automatically compact. It is not surprising that finite dimensionality
is crucial for this to happen.

\begin{theorem}
Let $\cl A$ and $\cl B$ be C*-algebras. The following statements are equivalent:

(i) \ Every element of $\frak{M}(\cl A,\cl B)$ is a compact multiplier;

(ii) Either $\cl A$ is finite dimensional and $\cl K(\cl B) = \cl B$
or $\cl B$ is finite dimensional and $\cl A = \cl K(\cl A)$.
\end{theorem}

The proof is based, in particular, on a result of Varopoulos
\cite{varopoulos} showing that if $X$, $Y$ are infinite compact
Hausdorff spaces then there exists a sequence
$(f_i)_{i=1}^{\infty}\subseteq C(X)\otimes_{\hh} C(Y)$ such that
$\sup_{i\in\bb{N}}\|f_i\|_{\hh} < \infty$ which converges uniformly
to a function $f\in C(X\times Y)\setminus C(X)\otimes_{\hh} C(Y)$.
By the Characterisation Theorem~\ref{charac}, such an $f$ must
belong to $\frak{M}(C(X),C(Y))$.

\section{Going multidimensional}

In the present section, we introduce a multidimensional version of
Schur and operator multipliers.

If $R_1,\dots,R_{n+1}$ are rings, $M_i$ is a $R_i$-left and
$R_{i+1}$-right module for each $i = 1,\dots,n$, and $M$ is an
$R_1$-left and $R_{n+1}$-right module, a multilinear map $\Phi :
M_1\times\dots \times M_n\rightarrow M$ is called
$R_1,\dots,R_{n+1}$-modular (or simply modular if
$R_1,\dots,R_{n+1}$ are clear from the context) if
$$\Phi(a_1 m_1 a_2,m_2 a_3,\dots,m_n a_{n+1})
= a_1\Phi(m_1, a_2 m_2, a_3 m_3,\dots,a_n m_n) a_{n+1},$$ for all
$m_i\in M_i$ ($i = 1,\dots,n$) and $a_j\in R_j$ ($j = 1,\dots,n+1$).

A multilinear Schur product was introduced by Effros and Ruan
\cite{effros_ruan_multi} as a multilinear map
$T:\underbrace{M_n(\mathbb C)\times\cdots \times M_n(\mathbb
C)}_n\to M_n(\mathbb C)$ which is
$\underbrace{D_n,\dots,D_n}_{n+1}$-modular, where $D_n$ is the
algebra of all diagonal matrices in $M_n(\mathbb C)$.

It is not difficult to see that any such mapping $T$ has the form
$$T(a^r,\ldots,a^1)_{i,j}=\sum_{(k_1,\ldots,k_{r-1})}A_{i,j}^{k_{r-1}\ldots k_1}a_{i,k_{r-1}}^r
a_{k_{r-1},k_{r-2}}^{r-1}\ldots a_{k_1j}^1,$$ where $a^i =
(a^i_{k,l})_{k,l}$ and, given the usual matrix units $e_{i,j}\in
M_n(\mathbb C)$,
$$A_{i,j}^{k_{r-1}\ldots k_1}=T(e_{i,k_{r-1}},
e_{k_{r-1},k_{r-2}}\ldots e_{k_1,j})_{i,j}.$$

The following theorem gives a characterisation of all bounded
multilinear Schur products.

\begin{theorem}
Suppose $T:\underbrace{M_n(\mathbb C)\times\cdots \times
M_n(\mathbb C)}_n\to M_n(\mathbb C)$ is a multilinear Schur product
map. Then the following are equivalent:

(i) \ the linearisation of $T$ is a contraction for the Haagerup
norm;

(ii) there exists a Hilbert space $H$, $2n$ contractions $a_1(j)\in
\cl B(H,\mathbb C)$,
 $a_r(i)\in \cl B(\mathbb C,H)$,
$i,j=1,\ldots,n$ and $n^{r-1}$ contractions $a_l(k)\in \cl B(H)$,
$l=2,\ldots,r-1$, $k=1,\ldots,n$ such that
$$A_{i,j}^{k_{r-1}\ldots k_1}=a_r(i)a_{r-1}(k_{r-1})\ldots a_2(k_1)a_1(j).$$
\end{theorem}

The theorem was proved in \cite{effros_ruan_multi} for complete
contraction $T$. A generalisation of Smith's result
\cite[Theorem~2.1]{smith} to the mutlidimensional setting,
\cite[Lemma~3.3]{jtt}, giving that any $D_n,D_n$-modular contraction $T$
is a complete contraction allows us to formulate the statement in
this generality. More about completely bounded multilinear maps will
be said later.

We now introduce multidimenisonal measurable Schur multipliers
following \cite{jtt}.

Let $(X_i,\mu_i)$, $i = 1,\dots,n$, be standard $\sigma$-finite
measure spaces. For notational convenience, integration with respect
to $\mu_i$ will be denoted by $dx_i$. Let
$$\Gamma(X_1,\dots,X_n) = L^2(X_1\times X_2)\odot L^2(X_2\times X_3)\odot\dots\odot L^2(X_{n-1}\times
X_n),$$  where each product $X_{i}\times X_{i+1}$ is equipped with the
corresponding product measure.

We identify the elements of $\Gamma(X_1,\dots,X_n)$ with functions
on $$X_1\times X_2\times X_2\times\dots\times X_{n-1}\times
X_{n-1}\times X_n$$ in the obvious fashion and  equip
$\Gamma(X_1,\dots,X_n)$ with two norms; one is the projective norm
$\|\cdot\|_{2,\wedge}$, where each of the $L^2$-spaces is equipped
with its $L^2$-norm, and the other is the Haagerup tensor norm
$\|\cdot\|_{\hh}$, where the $L^2$-spaces are given their opposite
operator space structure arising from the identification of
$L^2(X\times Y)$ with the class of Hilbert-Schmidt operators from
$L^2(X)$ into $L^2(Y)$ given by $f\mapsto T_f$.

For each $\nph\in L^{\infty}(X_1\times\dots\times X_n)$, we consider
a linear map $S_{\nph}$ defined on $\Gamma(X_1,\dots,X_n)$ and
taking values in $L^2(X_1\times X_n)$: for an
elementary tensor $f_1\otimes\dots\otimes f_{n-1}$ in
$\Gamma(X_1,\dots,X_n)$, we set $S_{\nph}(f_1\otimes\dots\otimes
f_{n-1})(x_1,x_n)$ to be equal to
\begin{equation}\label{multi}
\int_{X_2\times\dots\times X_{n-1}}
\nph(x_1,\dots,x_n)f_1(x_1,x_2)f_2(x_2,x_3)\dots
f_{n-1}(x_{n-1},x_n) dx_2\dots dx_{n-1}.
\end{equation}
One can show that $S_{\nph}$ is bounded as a map from
$(\Gamma(X_1,\dots,X_n),\|\cdot\|_{2,\wedge})$ into $(L^2(X_1\times
X_n),\|\cdot\|_2)$. Moreover, any multilinear bounded map
$S:L^2(X_1\times X_2)\times L^2(X_2\times X_3)\times\dots\times
L^2(X_{n-1}\times X_n)\rightarrow L^2(X_1\times X_n)$  which is
$L^\infty(X_1),\ldots L^\infty(X_n)$-modular  is given  by
(\ref{multi}) in analogy with Effros-Ruan's multilinear Schur product.

If, moreover, $S_\nph$ is bounded as a map from
$(\Gamma(X_1,\dots,X_n),\|\cdot\|_{\hh})$ into $(L^2(X_1\times
X_n),\|\cdot\|_{\op})$ that is, if there exists $C
> 0$ such that $\|S_{\nph}(F)\|_{\op}\leq C \|F\|_{\hh}$, for
all $F\in\Gamma(X_1,\dots,X_n)$, then we say that $\nph$ is a {\it
Schur $\mu_1,\dots,\mu_n$-multiplier} or simply a Schur multiplier,
if the measures are clear from the context. The smallest constant
$C$ with this property will be denoted by $\|\nph\|_{\mm}$.

We note also that if $H_i = L^2(X_i)$, $\cl D_i = \{M_{\psi} : \psi
\in L^{\infty}(X_i)\}$, $i = 1,\dots,n$, and
$$\hat{S}_{\nph} : \cl C_2(H_1,H_2)\times\dots\times \cl C_2(H_{n-1},H_n)
\rightarrow \cl C_2(H_1,H_n)$$ is the map defined by $\hat{S}_{\nph}
(T_{f_1},\dots,T_{f_{n-1}}) = T_{S_{\nph}(f_1,\dots,f_{n-1})}$, then if
$a_i\in L^{\infty}(X_i)$, $i = 1,\dots,n$, and $\nph(x_1,\ldots,
x_n)=a_1(x_1)\ldots a_n(x_n)$, we obtain
$$\hat{S}_{\nph} (T_{f_1},\dots,T_{f_{n-1}})=M_{a_n}T_{f_{n-1}}M_{a_{n-1}}\ldots T_{f_1}M_{a_1}.$$

Next theorem generalizes Theorem~\ref{groth2} to the
multidimensional case giving a characterisation of all Schur
multipliers in $L^\infty(X_1\times\ldots\times X_n)$.

\begin{theorem}\label{th_g3}
Let $\nph\in L^{\infty}(X_1\times\dots\times X_n)$. The following
are equivalent:

(i) \ $\nph$ is a Schur multiplier and $\|\nph\|_{\mm} < C$;

(ii) there exist essentially bounded functions $a_1 : X_1\rightarrow
M_{\infty,1}$, $a_n : X_n\rightarrow M_{1,\infty}$ and $a_i :
X_i\rightarrow M_{\infty}$, $i = 2,\dots,n-1$, such that, for almost
all $x_1,\dots,x_n$ we have
$$\nph(x_1,\dots,x_n) = a_n(x_n)a_{n-1}(x_{n-1})\dots a_1(x_1) \ \mbox{ and } \
\esssup_{x_i\in X_i}\prod_{i=1}^n \|a_i(x_i)\| < C.$$
\end{theorem}
The proof is based on the fact that if  $\nph$ is a Schur multiplier
then $\hat S_\nph$ gives rise to a multilinear map from $\cl
B(H_{n-1},H_n)\times\dots\times \cl B(H_1,H_2)$ into $\cl
B(H_1,H_n)$ which is completely bounded, normal and $\cl
D_n,\ldots,\cl D_1$-modular, and a characterisation of such maps
given by Christensen and Sinclair \cite{cs}.

\bigskip

The space of all functions satisfying condition (ii) of
Theorem~\ref{th_g3} can be identified with the extended Haagerup
tensor product
$L^\infty(X_1)\otimes_{\eh}\ldots\otimes_{\eh}L^\infty(X_n)$ which
will be discussed later.

In a same way two dimensional Schur multipliers are related to
double operator integrals, multidimensional Schur multipliers are
related to multiple operator integrals studied recently by Peller in
\cite{peller}. This notion is important due to its application to
the study of higher order differentiability of functions of
operators.

To define multiple operator integrals we fix spectral measures $
E_1(\cdot),\ldots,$ $E_n(\cdot)$ on $X_1,\ldots,X_n$, respectively.
Let $\mu_1,\ldots,\mu_n$ be scalar measures equivalent to
$E_1,\ldots, E_n$, respectively. Consider the space of all functions
$\nph$ for which there exists a measure space $(\mathcal T,\nu)$ and
measurable functions $g_i$ on $X_i\times\mathcal T$ such that
\begin{equation}\label{inproj}
\nph(x_1,\ldots,x_n)=\int_{\mathcal T}g_1(x_1,t)\ldots
g_n(x_n,t)d\nu(t),
\end{equation}
for almost all $x_1,\ldots,x_n$, where
$$\int_{\mathcal
T}\|g_1(\cdot,t)\|_{\infty}\ldots\|g_n(\cdot,t)\|_{\infty}d\nu(t)<\infty.$$
The space is called the integral projective tensor product of
$L^\infty(X_1),\ldots,$ $L^\infty(X_n)$ and denoted by
$L^\infty(X_1)\hat\otimes_i\ldots\hat\otimes_i L^\infty(X_n)$.

In the case $n=2$ this space coincides with the space of all Schur
multipliers by \cite{peller_two_dim}. For $n>2$ we can only show
that the space consists of  Schur multipliers.

For $\nph\in L^\infty(X_1)\hat\otimes_i\ldots\hat\otimes_i
L^\infty(X_n)$ and $(n-1)$-tuple of bounded operators
$(T_1,\ldots,T_{n-1})$ Peller defines a multiple operator integral
by
\begin{eqnarray*}
&&I_\nph((T_1,\dots,T_{n-1})=\\
&&\int_{\mathcal T}(\int_{X_1}g_1(x_1,t)d E_1(x_1))T_1
(\int_{X_2}g_2(x_2,t)d E_2(x_2))T_2\\
&&\ldots T_{n-1}(\int_{X_n}g_n(x_n,t) d E_n(x_n)))d\nu(t).
\end{eqnarray*}
If the spectral measures are multiplicity free and
$T_1,\ldots,T_{n-1}$ are Hilbert-Schmidt operators with respective
kernels $f_1,\ldots,f_{n-1}$ then one can easily see that
$I_\nph(T_1,\ldots, T_n)$ is a Hilbert-Schmidt operator with kernel
$S_\nph(f_1\otimes\ldots\otimes f_n)$.

\bigskip

Like measurable Schur multipliers, operator multipliers defined by
Kissin and Shulman can be generalised to the multidimensional
setting.

Let $H_1,\dots,H_n$ be Hilbert spaces ($n$ is even) and let $H =
H_1\otimes\dots\otimes H_n$. We define a Hilbert space
$HS(H_1,\dots,H_n)$ isometrically isomorphic to $H$: we let
$HS(H_1,H_2) = \cl C_2(H_1^{\dd},H_2)$, and by induction define
$$HS(H_1,\dots,H_n) = \cl C_2(HS(H_2,H_3)^{\dd},
HS(H_1,H_4,\dots,)).$$ Let $\theta: H\otimes K\to \cl
C_2(H^{\dd},K)$ be the natural surjective isometry from the product
of Hilbert spaces $H\otimes K$ to the Hilbert space $\cl
C_2(H^{\dd}, K)$ of Hilbert-Schmidt operators defined in
Section~\ref{sec3}. We extend this map by induction to the
multidimensional case to get a map $\theta : H \rightarrow
HS(H_1,\dots,H_n)$ by letting
$$\theta(\xi_{2,3}\otimes \xi) = \theta(\theta(\xi_{2,3})\otimes
\theta(\xi)),$$ where $\xi_{2,3} \in H_2\otimes H_3$ and $\xi\in
H_1\otimes H_4\otimes\dots\otimes H_n$.

Let $\Gamma(H_1,\dots,H_n) = (H_1\otimes H_2)\odot (H_2^{\dd}\otimes
H_3^{\dd})\odot\dots\odot (H_{n-1}\otimes H_n)$ equipped with the
Haagerup norm $\|\cdot\|_{\hh}$ where $H_i\otimes H_{i+1}$ is given
the operator space structure opposite to the one arising from the
embedding $\theta : H_i\otimes H_{i+1} \hookrightarrow \cl
B(H_i^{\dd},H_{i+1})$ (and similarly for $(H_i\otimes H_{i+1})^{\dd}
= H_i^{\dd}\otimes H_{i+1}^{\dd}$).

Fix $\nph\in B(H)$. We define a mapping $S_{\nph}:
\Gamma(H_1,\dots,H_n)\to \cl B(H_1^{\dd},H_n)$ as follows: if $\zeta
\in \Gamma(H_1,\dots,H_n)$ is an elementary tensor, namely,
$$\zeta =
\xi_{1,2}\otimes\eta_{2,3}^{\dd}\otimes\xi_{3,4}\otimes\dots\otimes\xi_{n-1,n},$$
we let
$$S_{\nph}(\zeta)=\theta(\nph(\xi_{1,2}\otimes\dots\otimes\xi_{n-1,n})(\theta(\eta_{2,3}^{\dd}))
\dots(\theta(\eta_{n-2,n-1}^{\dd})).$$

Using the natural identification, we consider $S_{\nph}$ as a map
from $\cl C_2(H_1^{\dd},H_2)\odot\ldots\odot \cl C_2(H_{n-1}^{\dd},
H_n)$ into $\cl B(H_1^{\dd},H_n)$ which in particular satisfies the
following:
$$S_{a_1\otimes\ldots\otimes a_n}(T_1\otimes\ldots\otimes
T_{n-1})=a_nT_{n-1}\ldots a_3^{\dd}T_2 a_2 T_1 a_1^{\dd}$$ for any
$a_1\otimes\ldots\otimes a_n\in \cl B(H)$.

For odd $n$, a similar definition can be given by \lq\lq adding'' to
$H_1,\dots,H_n$ the one-dimensional Hilbert space $\bb{C}$. For
technical simplicity from now on we restrict our attention to the case of even $n$
and refer the reader to \cite{jtt,jltt} for the general case.

We call $\nph$ a {\bf concrete operator multiplier} if  there exists
$C> 0$ such that
$$\|S_{\nph}(\zeta)\|_{\op}\leq C \|\zeta\|_{\hh}, \ \ \ \mbox{ for all } \zeta\in
\Gamma(H_1,\dots,H_n).$$

As in the two dimensional case, we want to specify classes of \lq\lq
continuous'' operator multipliers. Let $\cl A_i$ be a C*-algebra and
$\pi_i:\cl A_i\to B(H_i)$ be a representation, $i=1,\ldots,n$. An
element $\nph\in \cl A_1\otimes\ldots\otimes \cl A_n$ is  called
 an operator {\it $\pi_1,\dots,\pi_n$-multiplier} if
$(\pi_1\otimes\ldots\otimes\pi_n)(\nph)$ is a concrete operator
multiplier. We denote the set of all
$(\pi_1,\ldots,\pi_n)$-multipliers by $\frak{M}_{\pi_1,\dots,\pi_n}$
and denote by $\|\nph\|_{\pi_1,\ldots,\pi_n}$ the concrete
multiplier norm.

In analogy with the two-dimensional case we say that  $\nph$ is  a
{\it universal operator multiplier} if it is
$\pi_1,\dots,\pi_n$-multiplier for all choices of
$\pi_1,\dots,\pi_n$. In this case,
$$\|\nph\|_{\mm} \stackrel{def}{=} \sup
\|\nph\|_{\pi_1,\dots,\pi_n} < \infty.$$ By ${\mathfrak M}(\As)$ we
will denote the set of all universal operator multipliers. For $n =
2$ one obtains the notion of operator multipliers introduced by
Kissin and Shulman and discussed in Section~\ref{sec3}. Moreover,
natural analogs of the Comparision Theorem~\ref{compar} and
Corollary~\ref{cor_supp} hold in the new multidimensional setting.

As one may expect, multidimensional operator multipliers are related
to completely bounded  multilinear maps. We now recall the
definition and some related notions.

Let $\cl E$, $\cl E_1,\ldots,\cl E_n$ be closed subspaces of $\cl
B(H)$, $\cl B(H_1),\ldots, \cl B(H_n)$, respectively. We denote by
$\cl E_1\odot\ldots\odot\cl E_n$ the algebraic tensor product of
$\cl E_1,\ldots\cl E_n$. Let $a_k = (a_{i,j}^k)\in
M_{m_k,m_{k+1}}(\cl E_k)$, $k = 1,\dots,n$. We denote by
\begin{equation}\label{odot}
a^1\odot \dots\odot a^n\in M_{m_1,m_{n+1}}(\cl E_1\odot\dots\odot
\cl E_n)
\end{equation}
the matrix whose $i,j$-entry is
\begin{equation}\label{odot2}
\sum_{i_2,\dots,i_n} a^1_{i,i_2}\otimes
a^2_{i_2,i_3}\otimes\dots\otimes a^n_{i_n,j}.
\end{equation}
Let $\Phi : \cl E_1\times\dots\times \cl E_n \rightarrow\cl E$ be a
multilinear map and
$$\Phi^{(m)} : M_m(\cl E_1)\times M_m(\cl E_2)\times\dots\times M_m(\cl E_n)
\rightarrow M_m(\cl E)$$ be the multilinear map given by
\begin{equation}\label{rfl}
\Phi^{(m)}(a^1,\dots,a^n)_{i,j} = \sum_{i_2,\dots,i_n}
\Phi(a^1_{i,i_2},a^2_{i_2,i_3},\dots, a^n_{i_n,j}),
\end{equation}
where $a^k = (a^k_{i,j})\in M_m(\cl E_k)$, $1\leq i,j\leq m$. The
map $\Phi$ is called completely bounded if there exists $C
> 0$ such that for all $m\in\bb{N}$ and all elements $a^k\in
M_{m}(\cl E_k)$, $k = 1,\dots,n$, we have
$$\|\Phi^{(m)}(a^1,\dots,a^n)\|\leq C \|a^1\|\dots\|a^n\|.$$

Every completely bounded multilinear map $\Phi : \cl
E_1\times\dots\times \cl E_n \rightarrow\cl E$ gives rise to a
completely bounded linear map from the Haagerup tensor product $\cl
E_1\otimes_{\hh}\dots\otimes_{\hh} \cl E_n$ into $\cl E$.

The {\it extended Haagerup tensor product} ${\mathcal
E}_1\ehaags{\mathcal E}_n$ is defined in \cite{effros_ruan} as the
space of all normal (in each variable) completely bounded maps
$u:{\cl E_1}^*\timess{\cl E_n}^*\to{\mathbb C}$. It was shown
in~\cite{effros_ruan} that if $u\in {\cl E_1}\ehaags{\cl E_n}$ then
there exist index sets $J_1,J_2,\dots,J_{n-1}$ and matrices $a^1 =
(a^1_{1,s})\in M_{1,J_1}({\cl E_1})$, $a^2 = (a^2_{s,t})\in
M_{J_1,J_2}({\cl E_2}),\dots, a^n = (a^n_{t,1})\in
M_{J_{n-1},1}({\cl E_n})$ such that if $f_i\in{\cl E_i}^*$,
$i=1,\dots,n$, then
\begin{equation}\label{product}
\langle u,f_1\otimess f_n\rangle\defeq u(f_1,\dots,f_n) = \langle
a^1,f_1\rangle\dots\langle a^n,f_n\rangle,
\end{equation}
where $\langle a^k,f_k\rangle=\big(f_k(a_{s,t}^k)\big)_{s,t}$ and
the product of the (possibly infinite) matrices in~\eqref{product}
is defined to be the limit of the sums
$$\sum_{i_1\in F_1,\dots, i_{n-1}\in F_{n-1}}f_1(
a_{1,i_1}^1)f_2(a^2_{i_1,i_2})\dots f_n(a_{i_{n-1},1}^n)$$ along the
net $\{(F_1\timess F_{n-1}) : F_j\subseteq J_j \mbox{ finite}, 1\leq
j\leq n-1\}$.

We identify $u$ with the matrix product $a^1\odots a^n$; two
elements $a^1\odots a^n$ and $\tilde a^1\odots \tilde a^n$ coincide
if $\langle a^1,f_1\rangle\dots\langle a^n,f_n\rangle = \langle
\tilde a^1,f_1\rangle\dots\langle \tilde a^n,f_n\rangle$ for all
$f_i\in{\cl E_i}^*$, $i = 1,\dots,n$. Moreover,
$$\|u\|_\eh = \inf\{\|a^1\|\dots\|a^n\| : u = a^1\odots
a^n\}.$$

There is a natural bijection $\gamma$ between the extended Haagerup
tensor product $\cl B(H_1)\ehaags \cl B(H_n)$ and the space of
multilinear normal completely bounded maps from $\cl B(H_2,
H_1)\times\ldots\times \cl B(H_n, H_{n-1})$ to $\cl B(H_n, H_1)$
given as follows: if $u = A_1\odot \dots\odot A_n \in \cl
B(H_1)\ehaags \cl B(H_n)$ then
$$\gamma(u)(T_1,\dots,T_{n-1}) = A_1 (T_1\otimes I)A_2\dots
A_{n-1}(T_{n-1}\otimes I)A_n,$$ for all $T_i\in \cl B(H_{i+1},H_i)$,
$i = 1,\dots,n-1$. This is due to  Christensen and
Sinclair~\cite{cs}.

The connection of multilinear completely bounded maps with universal
multidimensional operator multipliers arises as follows. Let $\cl
A_i$ be a $C^*$-algebra, $i=1,\ldots,n$, and $\nph\in {\mathfrak
M}(\cl A_1,\ldots, \cl A_n)$. Then the map $S_{\nph}$ is completely
bounded for the opposite operator space structures, and hence has a
completely bounded extension to a map $$\Phi_\nph:\left(\cl
K(H_{n-1}^{\dd},H_n)\haags\cl
K(H_1^{\dd},H_2),\|\cdot\|_{\hh}\right)\rightarrow \left(\cl
K(H_1^{\dd},H_n),\|\cdot\|_{\op}\right)$$ given by
$$\Phi_{\nph}(T_{n-1}\otimes\ldots\otimes T_1)=S_{\nph}(T_1\otimes\ldots\otimes T_{n-1}).$$
Thus, $\Phi_\nph^{**}$ is a completely bounded normal map, and the
$(\cl A_n', (\cl A_{n-1}^{\dd})',\ldots,$ $(\cl
A_1^{\dd})')$-modularity of $\Phi_\nph$ allows then to define a
symbol of a multidimensional universal multiplier:

\begin{theorem}\label{symb}
Let $\As$ be C*-algebras and $\nph\in M(\As)$. There exists a unique
element $u_{\nph}\in
\A_n\ehaag\A_{n-1}^{o}\ehaags\A_2\ehaag\A_1^{o}$ with the property
that if $\pi_i$ is a representation of $\A_i$ for $i = 1,\dots,n$
then the map $\Phi_{\pi_1\dots\pi_n(\nph)}$ coincides with the
restriction of
$\gamma(\pi_n\otimes_{\eh}\pi_{n-1}^{\dd}\otimes_{\eh}\ldots\otimes_{\eh}\pi_1^{\dd}(u_{\nph}))$.

The map $\nph\rightarrow u_{\nph}$ is linear and if $a_i\in\A_i$, $i
= 1,\dots,n$ then $u_{a_1\otimess a_n} = a_n\otimes
a_{n-1}^o\otimess a_1^o.$ Moreover, $\|\nph\|_{\mm} =
\|u_{\nph}\|_{\eh}$.
\end{theorem}

The notion of completely compact map, completely compact  and
compact multipliers have natural extensions to the mutlidimensional
case. Namely,  if $\Y, \X_1,\dots,\X_n$ are operator spaces and
$\Phi : \X_1\timess\X_n\rightarrow\Y$ is a completely bounded
multilinear map, we call $\Phi$ {\it completely compact} if for each
$\epsilon
> 0$ there exists a finite dimensional subspace $F\subseteq\Y$
such that
\[\dist(\Phi^{(m)}(x_1,\dots,x_n), M_m(F)) < \epsilon,\] for all
$x_i\in M_{m}(\X_i)$, $\|x_i\|\leq 1$, $i = 1,\dots,n$, and all
$m\in\bb{N}$. We
denote by $\CC(\X_1\timess\X_n, \Y)$ %
the space of all completely bounded completely compact multilinear
maps from $\X_1\timess\X_n$ into $\Y$.

Let $\A_i$ be a $C^*$-algebra, $i=1,\ldots,n$. We define
$\frak{M}_{cc}(\As)$ (resp. $\frak{M}_{c}(\As)$ and
$\frak{M}_{\ff}(\As)$) as the set of all $\nph\in {\mathfrak
M}(\As)$ such that there exist faithful representations
$\pi_1,\ldots,\pi_n$ of $\A_1,\ldots,\A_n$ with the property that if
$\pi = \pi_1\otimes\ldots\pi_n$ then $\Phi_{\pi(\nph)}$ is
completely compact (resp. $\Phi_{\pi(\nph)}$ is compact and the
range of $\Phi_{\pi(\nph)}$ is a finite dimensional space of
finite-rank operators).

Let \[\Kh\defeq\K(H_2,H_1)\haags \K(H_{n},H_{n-1}).\] Saar's result
(Theorem~\ref{saarcc}) has the following generalisation:
\begin{theorem}\label{muccc}
   The operator space
  $\E\defeq \K(H_1)\haag(\B(H_2)\ehaags\B(H_{n-1}))\haag \K(H_n)$ is isometrically isomorphic to
  $\F\defeq \CC(\Kh,\K(H_n,H_1))$ with an isometry given by the restriction of $\gamma$ to $\E$.
\end{theorem}

This leads to the following characterisation of completely compact
universal multipliers in terms of their symbols:

\begin{theorem}\label{charcn}
  Let $\As$ be C*-algebras
  and $\nph\in M(\As)$. The following are
  equivalent:\smallskip

(i) \ \ $\nph\in M_{cc}(\As)$;\smallskip

(ii) \ $u_{\nph}\in \K(\A_n) \haag(\A_{n-1}^o\ehaags\A_2)\haag
\K(\A_1^o)$

(iii) there exists a net $\{\nph_{\alpha}\}_{\alpha}\subseteq
M_{\ff}(\As)$ such that $\|\nph_{\alpha} -
\nph\|_{\mm}\rightarrow_{\alpha} 0$.
\end{theorem}

The important point to note here is that there is no direct analogy
with the two-dimensional case: the space $\CC(\Kh,\K(H_n,H_1))$ is
not isometrically isomorphic to $\K(H_1)\haag\K(H_2)\haag\ldots\haag
\K(H_n)$, and the  symbol of a completely compact map maybe an
element of a bigger  than $\K(\A_n)
\haag\K(\A_{n-1}^o)\haag\ldots\haag \K(\A_1^o)$ space contrary to
what one may expect by following the analogy with the
two-dimensional case.

The property  of the set of  universal completely compact
multipliers in $\cl A_1\otimess\cl A_n$  to coincide with  the set
of compact universal multipliers depends only on the structure of
compact elements of $\A_1$ and $\A_n$  and not on $\A_k$,
$k=2,\ldots,n-1$. The condition on $\K (\A_1)$ and $\K (\A_n)$ are
exactly the ones given in Theorem~\ref{th_automcc}.

We shall end the section by describing  an application of
multidimensional Schur multipliers to abstract harmonic analysis.
Let $G$ be a locally compact, $\sigma$-compact group and let $A(G)$
be  the Fourier algebra of $G$. We recall that $\lambda_s$  denote
the left regular representation of $G$ on $L^2(G)$.  Since $A(G)$ is
the predual of the von Neumann algebra $\vn(G)$, it possesses a
canonical operator structure. Therefore we can define the
multidimensional Fourier algebra as follows:
$$A^n(G) = \underbrace{A(G)\ehaags A(G)}_n.$$
Since $\vn(G)$ is generated by $\lambda_x$, $x\in G$, it follows
from the definition that the elements of $A^n(G)$ can be identified
with functions $f\in L^\infty(G^n)$ given by
$f(x_n,\ldots,x_1)=\Phi(\lambda_{x_n},\ldots,\lambda_{x_1})$, where
$\Phi : \vn(G)\timess \vn(G)\to \mathbb C$ is a normal completely
bounded multilinear map.

For $f\in A(G)$ we define a map $\theta:A(G)\to A^n(G)$ by
$$\theta(f)(x_n,\ldots,x_1)=f(x_n\ldots x_1).$$
We call a function $\nph\in L^\infty(G^n)$ an {\it $n$-multiplier}
of $A(G)$ if $\nph\theta(f)\in A^n(G)$ whenever $f\in A(G)$. We call
$\nph$  {\it a completely bounded $n$-multiplier} if the mapping
$f\mapsto\nph\theta(f) $ is completely bounded. The following
characterisation, which is a multidimensional version of Theorem~\ref{bf},
was established in \cite{tt}.

\begin{theorem}
Let $\nph\in L^\infty(G^n)$. The following are equivalent:

(i) \ $\nph$ is a completely bounded $n$-multiplier;

(ii) The function $\tilde\nph\in L^\infty(G^{n+1})$ given by
$$\tilde\nph(x_1,\ldots,x_n)=\nph(x_{n+1}^{-1}x_n,\ldots,x_2^{-1}x_1)$$
is a Schur multiplier with respect to the left Haar measure.
\end{theorem}

\end{document}